\documentclass[12pt,a4paper]{amsart}
\usepackage{amssymb,latexsym,graphicx}
\usepackage{color}
\usepackage[cp1252]{inputenc} 
\usepackage[T1]{fontenc}
\usepackage{lmodern}
\usepackage[pagebackref=true]{hyperref}
\hypersetup{pdftex,colorlinks=true,linkcolor=blue,citecolor=blue,urlcolor=blue}
\usepackage{currfile}
\usepackage[english]{babel}
\usepackage[section]{placeins}
\usepackage{version}

\excludeversion{hide}

\theoremstyle{plain} 

\newtheorem{theorem}{Theorem}[section]
\newtheorem{thm}[theorem]{Theorem}
\newtheorem{proposition}[theorem]{Proposition}
\newtheorem{corollary}[theorem]{Corollary}
\newtheorem{lemma}[theorem]{Lemma}

\newtheorem{conjecture}[theorem]{Conjecture}
\newtheorem{questions}[theorem]{Questions}

\theoremstyle{definition} 

\newtheorem{remark}[theorem]{Remark}
\newtheorem{remarks}[theorem]{Remarks}

\numberwithin{equation}{section}
\numberwithin{figure}{section}
\numberwithin{table}{section}

\newcommand{\noib}{\noindent $\diamond$~} 

\definecolor{purple}{RGB}{127,0,255}

\newcommand{\D}{\mathbb{D}}

\newcommand{\N}{\mathbb{N}}
\newcommand{\Nb}{\mathbb{N^{\bullet}}}

\newcommand{\R}{\mathbb{R}}
\newcommand{\Ss}{\mathbb{S}}
\newcommand{\bS}{\mathbb{S}}
\newcommand{\T}{\mathbb{T}}
\newcommand{\Z}{\mathbb{Z}}

\newcommand{\cD}{\mathcal{D}}
\newcommand{\cE}{\mathcal{E}}

\newcommand{\cK}{\mathcal{K}}
\newcommand{\cL}{\mathcal{L}}

\newcommand{\cP}{\mathcal{P}}

\newcommand{\cS}{\mathcal{S}}
\newcommand{\cT}{\mathcal{T}}

\newcommand{\cZ}{\mathcal{Z}}

\newcommand{\ip}[2]{\langle #1|#2 \rangle}

\newcommand{\sm}{\!\setminus\!}

\newcommand{\set}[1]{\left\{ #1 \right\}}
\DeclareMathOperator{\supp}{Supp}
\DeclareMathOperator{\sign}{sign}

\begin{document}


\title[]{Non-boundedness of the number of super level domains of eigenfunctions}

\author[P. B\'{e}rard]{Pierre B\'erard}
\author[P. Charron]{Philippe Charron}
\author[B. Helffer]{Bernard Helffer}

\address{PB: Universit\'{e} Grenoble Alpes and CNRS\\
Institut Fourier, CS 40700\\ 38058 Grenoble cedex 9, France.}
\email{pierrehberard@gmail.com}

\address{PC: Universit\'e de Montr\'eal\,\\
2920, Chemin de la Tour,
Montr\'eal, QC,
H3T 1J4, Canada}
\email{charronp@dms.umontreal.ca}

\address{BH: Laboratoire Jean Leray, Universit\'{e} de Nantes and CNRS\\
F44322 Nantes Cedex, France and LMO (Universit\'e Paris-Sud).}
\email{Bernard.Helffer@univ-nantes.fr}

\thanks{The authors would like to thank T.~Hoffmann-Ostenhof and E.~Lieb for motivating discussions on the subject, and I.~Polterovich for his comments on an earlier version. They are very grateful to the anonymous referee for her/his careful reading, and for the many comments which helped them improve this paper.}

\keywords{Eigenfunction, Nodal domain, Courant nodal domain theorem.}

\subjclass[2010]{35P99, 35Q99, 58J50.}

\date{\today ~(\currfilename)}

\begin{abstract} Generalizing Courant's nodal domain theorem, the ``Extended Courant property'' is the statement that a linear combination of the first $n$ eigenfunctions has at most $n$ nodal domains. A related question is to estimate the number of connected components of the (super) level sets of a Neumann eigenfunction $u$. Indeed, in this case, the first eigenfunction is constant, and looking at the level sets of $u$ amounts to looking at the nodal sets $\set{u-a=0}$, where $a$ is a real constant. In the first part of the paper, we prove that the Extended Courant property is false for the subequilateral triangle and for regular $N$-gons ($N$ large), with the Neumann boundary condition. More precisely, we prove that there exists a Neumann eigenfunction $u_k$ of the $N$-gon, with labelling $k$, $4 \le k \le 6$, such that the set $\set{u_k \not = 1}$ has $(N+1)$ connected components.
In the second part, we prove that there exists a metric $g$ on $\mathbb{T}^2$ (resp. on $\mathbb{S}^2$), which can be chosen arbitrarily close to the flat metric (resp. round metric), and an eigenfunction $u$ of the associated Laplace-Beltrami operator, such that the set $\set{u \not = 1}$ has infinitely many connected components. In particular the Extended Courant property is false for these closed surfaces. These results are strongly motivated by a recent paper by Buhovsky, Logunov and Sodin.  As for the positive direction, in Appendix~B, we prove that the Extended Courant property is true for the isotropic quantum harmonic oscillator in $\mathbb{R}^2$.
\end{abstract}%

\maketitle

\section{Introduction}\label{S-int}

Let $\Omega \subset \R^2$ be a bounded domain (open connected set) with piecewise smooth boundary,  or a compact Riemannian surface, with or without boundary, and let $\Delta$ be the  Laplace-Beltrami operator. Consider the (real) eigenvalue problem
\begin{equation}\label{E-int-2}
\left\{
\begin{array}{l}
- \Delta u = \lambda u \text{~in~} \Omega\,,\\[5pt]
B(u) = 0 \text{~on~} \partial \Omega\,,
\end{array}%
\right.
\end{equation}
where the boundary condition $B(u)=0$ is either the Dirichlet or the Neumann boundary condition, $u=0$ or $\partial_{\nu}u=0$ on $\partial \Omega\,$,  or the empty condition if $\partial \Omega$ is empty. \medskip

We arrange the eigenvalues of \eqref{E-int-2} in nondecreasing order, multiplicities taken into account,
\begin{equation}\label{E-in-2a}
\lambda_1 < \lambda_2 \le \lambda_3 \le \cdots
\end{equation}

The \emph{nodal set} $\cZ(u)$ of a (real) function $u$ is defined to be
\begin{equation}\label{E-int-4}
\cZ(u) = \overline{\set{x \in \Omega ~|~ u(x)=0}}\,.
\end{equation}

The \emph{nodal domains} of a function $u$ are the connected components of $\Omega \sm \cZ(u)$. Call $\beta_0(u)$ the number of nodal domains of the function $u$.\medskip

The following classical theorem can be found in \cite[Chap. VI.6]{CH1953}.

\begin{thm}[Courant, 1923]\label{T-CT}
An eigenfunction $u$, associated with the $n$-th eigenvalue $\lambda_n$ of the eigenvalue problem \eqref{E-int-2}, has at most $n$ nodal domains, $\beta_0(u) \le n$.
\end{thm}%

For $n \ge 1$, denote by $\cL_n(\Omega)$ the vector space of linear combinations of eigenfunctions of problem \eqref{E-int-2}, associated with the $n$ first eigenvalues, $\lambda_1, \dots, \lambda_n$.

\begin{conjecture}[Extended Courant Property]\label{C-ECP}
Let $w \in \cL_n(\Omega)$ be a nontrivial linear combination of eigenfunctions associated with the $n$ first eigenvalues of problem \eqref{E-int-2}. Then, $\beta_0(w) \le n$.
\end{conjecture}%

This conjecture is motivated by a statement made in a footnote\footnote{p. 454 in \cite{CH1953}.} of  Courant-Hilbert's book.\medskip

Conjecture~\ref{C-ECP} is known to be true in dimension $1$ (Sturm, 1833).
In higher dimensions, it was pointed out by V.~Arnold (1973), in relation with Hilbert's 16th problem, see \cite{Arn2014}. Arnold noted that the conjecture is true for $\R\mathrm{P}^2$, the real projective space with the standard metric. It follows from \cite{Ley1996} that Conjecture~\ref{C-ECP} is true when restricted to linear combinations of even (resp. odd) spherical harmonics on $\Ss^2$ equipped with the standard metric. Counterexamples  to the conjecture  were constructed by O.~Viro (1979) for  $\R\mathrm{P}^3$, see \cite{Vir1979}. As far as we know, $\R\mathrm{P}^2$ is the only higher dimensional compact  example for which Conjecture~\ref{C-ECP} is proven to be true. In Appendix~\ref{A-QHO}, we prove that the conjecture is true for the isotropic quantum harmonic oscillator in $\R^2$ as well. Simple counterexamples to Conjecture~\ref{C-ECP} are given in \cite{BH2018-ecp1,BH-teqa,BH-ecp-bal}. They include smooth convex domains in $\R^2$, with Dirichlet or Neumann boundary conditions.  A question related to the Extended Courant property is to estimate the number of connected components of the (super) level sets of a Neumann eigenfunction $u$. Indeed, in this case, the first eigenfunction is constant, and looking at the level sets of $u$ amounts to looking at the nodal sets $\set{u-a=0}$, where $a$ is a real constant. Most counterexamples to the Extended Courant property, not all, are of this type. This is the case in the present paper\footnote{We changed the initial title of our paper (arXiv:1906.03668v2, June 20, 2019) to reflect this fact, as suggested by the referee.} as well. Studying the topology of level sets of a Neumann eigenfunction is, in itself, an interesting question which is related to the ``hot spots'' conjecture, see \cite{BaPa2006}.

\begin{questions}\label{Q-ECP}
Natural questions, related to Conjecture~\ref{C-ECP}, arise.
\begin{enumerate}
  \item Fix $\Omega$ as above, and $N \ge 2$. Can one bound $\beta_0(w)$, for $w \in \cL_{N}(\Omega)$, in terms of $N$ and geometric invariants of $\Omega$?
  \item Assume that $\Omega \subset \R^2$ is  a convex domain. Can one bound $\beta_0(w)$, for $w \in \cL_N(\Omega)$, in terms of $N$,  independently of $\Omega$?
  \item Assume that $\Omega$ is a simply-connected closed surface. Can one bound $\beta_0(w)$, for $w \in \cL_N(\Omega)$, in terms of $N$,  independently of $\Omega$?
\end{enumerate}
\end{questions}%

 A negative answer to Question~\ref{Q-ECP}(1) for the $2$-torus is given in \cite{BLS}. In that paper, the authors construct a smooth metric $g$ on $\mathbb{T}^2$, and a family of eigenfunctions $\phi_j$ with infinitely many isolated critical points. As a by-product of their construction, they prove that there exist a smooth metric $g$, a family of eigenfunctions $\phi_j$, and a family of real numbers $c_j$ such that  $\beta_0(\phi_j-c_j)= +\infty$, see Proposition~\ref{P-BLS}.\medskip


The main results of the present paper are as follows. In Section~\ref{S-subteq}, we prove that Conjecture~\ref{C-ECP} is false for a subequilateral (to be defined later on) triangle with Neumann boundary condition, see Proposition~\ref{P-isos-2}. In Section~\ref{S-polyN}, we prove that the regular $N$-gons, with Neumann boundary condition, provide negative answers to both Conjecture~\ref{C-ECP}, and  Question~\ref{Q-ECP}(2),  at least for $N$ large enough, see Proposition~\ref{P-isos-4}. \medskip

 The second part of the paper, Sections~\ref{S-T2} and \ref{S-S2m}, is  strongly motivated by \cite{JaNa1999,BLS}. We give a new proof that Conjecture~\ref{C-ECP} is false for the torus $\T^2$, and we prove that it is false for the sphere $\Ss^2$ as well.  More precisely, we prove the existence of a smooth metric $g$ on $\T^2$ (resp. $\Ss^2$), which can be chosen arbitrarily close to the flat metric (resp. round metric), and an eigenfunction $\Phi$ of the associated Laplace-Beltrami operator, such that the set $\set{\Phi \not =1}$ has infinitely many connected components. We refer to Proposition~\ref{P-ex1-2} for the torus, and to Propositions~\ref{P-S2m-2} and \ref{P-S2m-8} for the sphere. \medskip

In the case of $\T^2$, we also consider real \emph{analytic} metrics. For such a metric, an eigenfunction can only have finitely many isolated critical points. In \cite[Introduction]{BLS}, the authors ask whether, for analytic metrics, there exists an asymptotic upper bound for the number of critical points of an eigenfunction, in terms of the corresponding eigenvalue. Proposition~\ref{P-ex2-2} is related to this question. In Section~\ref{S-fin}, we make some final comments.

\medskip

In Appendix~\ref{A-S2}, we prove the weaker result $\beta_0(w) \le 8\, d^2$ when $w$ is the restriction to $\Ss^2$ of a polynomial of degree $d$ in $\R^3$. This gives a partial answer to Question~\ref{Q-ECP}(3) in the case of the sphere. In Appendix~\ref{A-QHO}, we prove that Conjecture~\ref{C-ECP} is true for the isotropic quantum harmonic oscillator in $\R^2$. Both appendices rely on  \cite{Char2015}.

\section{Subequilateral triangle, Neumann boundary condition}\label{S-subteq}

\begin{figure}[!hbt]
\centering
\includegraphics[scale=0.4]{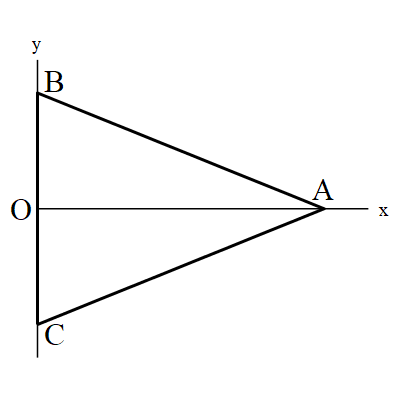}\vspace{-6mm}
\caption{Subequilateral triangle, $BC < AB=AC$}\label{F1}
\end{figure}

Let $\cT(b)$ denote the interior of the triangle with vertices $A=(\sqrt{3},0)$, $B=(0,b)$, and $C=(0,-b)$. When $b=1$, $\cT(1)$ is an equilateral triangle with sides of length $2$. From now on, we assume that $0 < b < 1$. The angle at the vertex $A$ is less than $\frac{\pi}{3}$, and we say that $T(b)$ is a \emph{subequilateral} triangle, see Figure~\ref{F1}. Let $\cT(b)_+ = \cT(b) \cap \{y > 0\}$, and $\cT(b)_- = \cT(b) \cap \{y < 0\}$.

\FloatBarrier

Call $\nu_i(\cT(b))$ the Neumann eigenvalues of $\cT(b)$, and write them in non-decreasing order, with multiplicities, starting from the labelling $1$,
\begin{equation}\label{E-iso-2}
0 = \nu_1(\cT(b)) < \nu_2(\cT(b)) \le \nu_3(\cT(b)) \le \cdots
\end{equation}

We recall the following theorems.

\begin{theorem}[\cite{LaSi2010}, Theorem~3.1]\label{T-isos-LS}
Every second Neumann eigenfunction of a subequilateral triangle $\cT(b)$ is even in $y$, $u(x,-y) = u(x,y)$.
\end{theorem}%

\begin{theorem}[\cite{Miya2013}, Theorem~B]\label{T-isos-M}
Let $\cT(b)$ be a subequilateral triangle. Then, the eigenvalue $\nu_2(\cT(b))$ is simple, and an associated eigenfunction $u$ satisfies $u(O) \not = 0$, where $O$ is the point $O=(0,0)$. Normalize $u$ by assuming that $u(O) = 1$. Then, the following properties hold.
\begin{enumerate}
  \item The partial derivative $u_x$ is negative in $\overline{\cT(b)}\setminus \left( \overline{BC} \cup \{A\} \right)$.
  \item The partial derivative $u_y$ is positive in $\overline{\cT(b)_+}\setminus \left( \overline{OA} \cup \{B\}\right)$, and negative in
      $\overline{\cT(b)_-}\setminus \left( \overline{OA} \cup \{C\}\right)$.
  \item The function $u$ has exactly four critical points $O,A,B$ and $C$ in $\overline{\cT}$.
  \item The points $B$ and $C$ are the global maxima of $u$, and $u(B)=u(C) > u(O) > 0$.
  \item The point $A$ is the  global minimum of $u$, and $u(A) < 0$.
  \item The point $O$ is the saddle point of $u$.
\end{enumerate}
\end{theorem}%

As a direct corollary of these theorems, we obtain the following result.

\begin{figure}[!hbt]
\centering
\includegraphics[scale=0.3]{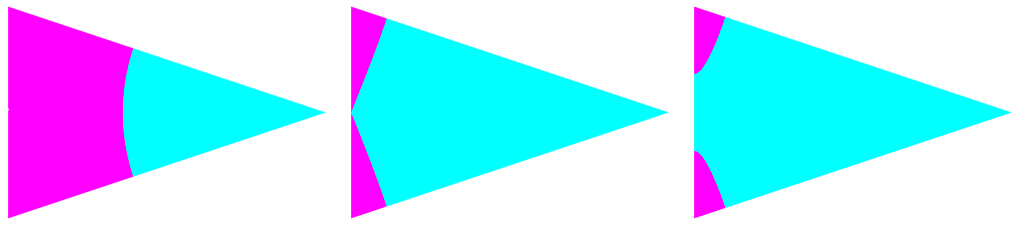}
\caption{Nodal domains of $u-a$ ($a<1,\, a=1,\, 1<a$)}\label{F2a}
\end{figure}

\begin{proposition}\label{P-isos-2}
Let $u$ be the second Neumann eigenfunction of the subequilateral triangle $\cT(b)$, $0 < b < 1$, normalized so that,
\begin{equation*}
u(A) = \min u < 0 < u(O) = 1 < \max u = u(B)=u(C) \,.
\end{equation*}
For $a \in \R$, let $\beta_0(u-a)$ be the number of nodal domains of the function $u-a$ (equivalently the number of $a$-level domains of $u$). Then,
\begin{equation*}
\left\{
\begin{array}{lll}
\beta_0(u-a) = 1 & \text{for} & a \le \min u \text{~or~} a \ge \max u\,,\\[5pt]
\beta_0(u-a) = 2 & \text{for} & \min u <  a < 1\,,\\[5pt]
\beta_0(u-a) = 3 & \text{for} &  1 \le a  < \max u\,.
\end{array}
\right.
\end{equation*}
As a consequence, for $1 \le a  < \max u$, the linear combination $u-a$ provides a counterexample to Conjecture~\ref{C-ECP}, see Figure~\ref{F2a}.
\end{proposition}%

\proof Fix some $0 < b < 1$, denote $\cT(b)$ by $\cT$, and $\nu_2(\cT(b))$ by $\nu_2$  for simplicity. In the proof, we write $(A1)$ for Assertion (1) of Theorem~\ref{T-isos-M}, etc..\medskip

For $a \in \R$, call $v_a$ the function $v_a := u - a$. This is a linear combination of a second and first Neumann eigenfunctions of $\cT$. We shall now describe the nodal set of $v_a$ carefully.\medskip

\begin{figure}[!hbt]
\centering
\includegraphics[scale=0.25]{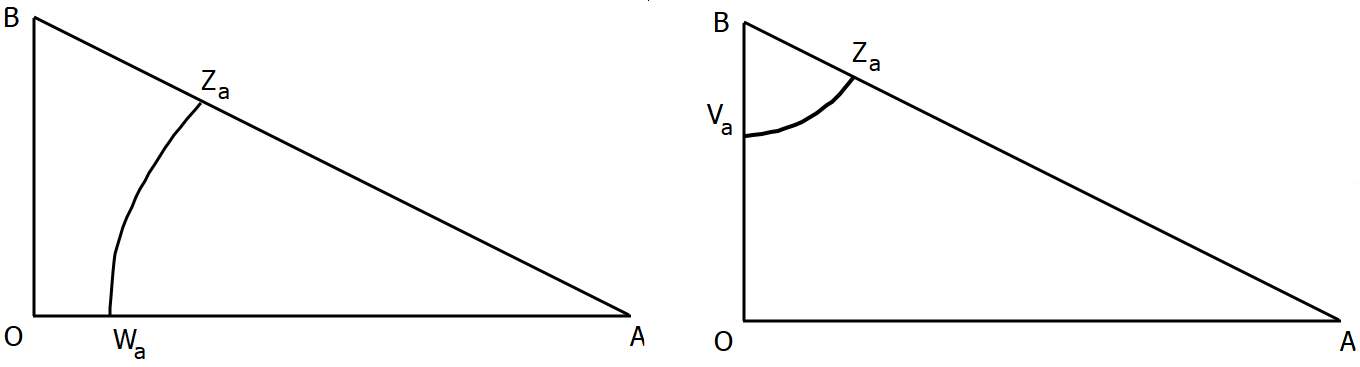}
\caption{Nodal behaviour of $u-a$ in the triangle $OAB$}\label{F2}
\end{figure}

According to Theorem~\ref{T-isos-LS}, for all $a$, the function $v_a$ is even in $y$, so that it is sufficient to determine its nodal set in the triangle $\cT_+ = OAB$, see Figure~\ref{F2}.\medskip

\noib By (A4) and (A5), the nodal set $\cZ(v_a)$ is nontrivial if and only if $u(A) < a < u(B)$.\smallskip

\noib By (A1) and (A2), the directional derivative of $v_a$ in the direction of $\overset{\longrightarrow}{BA}$ is negative in the open segment BA, so that $v_a|_{BA}$ is strictly decreasing from $v_a(B)$ to $v_a(A)$, and therefore vanishes at a unique point $Z_a = (\xi_a,\eta_a) \in BA$.  We now consider three cases.\smallskip

\textbf{Case $u(A) < a < u(O)$.}\smallskip

\noib By (A1), $v_a|_{OA}$ is strictly decreasing from $v_a(O)$ to $v_a(A)$, and therefore vanishes at a unique point $W_a = (\omega_a,0) \in OA$. By (A2), $\omega_a < \xi_a$.\smallskip

\noib By (A1) and (A2), it follows that the nodal set $\cZ(v_a)\cap \cT_{+}$ is contained in the rectangle $[\omega_a,\xi_a]\times [0,\eta_a]$, and that it is a smooth $y$-graph over $[\omega_a,\xi_a]$, and a smooth $x$-graph over $[0,\eta_a]$.\smallskip

We have proved that $v_a$ has exactly two nodal domains in $\cT$.\smallskip

\textbf{Case $a_c = u(O)$.}\smallskip

The analysis is similar to the previous one, except that $\omega_{a_c}=0$. As a consequence, $v_{a_c}$ has exactly three nodal domains in $\cT$.\smallskip

\textbf{Case $u(O) < a < u(B)$.}\smallskip

\noib By (A2), $v_a|_{OB}$ is strictly increasing from $v_a(0)$ to $v_a(B)$, so that it vanishes at a unique point $V_a = (0,\zeta_a) \in OB$. From (A1), it follows that $\zeta_a < \eta_a$. \smallskip

\noib From (A1) and (A2), it follows that the nodal set $\cZ(v_a)\cap \cT_{+}$ is contained in the rectangle $[0,\xi_a]\times [\zeta_a,\eta_a]$, and that it is a smooth $y$-graph over $[0,\xi_a]$, and a smooth $x$-graph over $[\zeta_a,\eta_a]$.\smallskip

It follows that $v_a < 0$ in $]-\zeta_a,\zeta_a[ \times [0,\sqrt{3}]\cap \cT$, and that $v_a$ has precisely three nodal domains in $\cT$.  Proposition~\ref{P-isos-2} is proved. \hfill\qed \medskip

\FloatBarrier

\section{Regular $N$-gon, Neumann boundary condition}\label{S-polyN}

\begin{proposition}\label{P-isos-4}
Let $\cP_N$ denote the regular polygon with $N$ sides, inscribed in the unit disk $\cD$. Then, for $N$ large enough, Conjecture~\ref{C-ECP} is false for $\cP_N$, with the Neumann boundary condition. More precisely, there exist $m \le 6$, an eigenfunction $u_m$ associated with $\nu_m(\cP_N)$, and a value $a$ such that the function $u_m - a$ has $N+1$ nodal domains.
\end{proposition}%

\proof  The general idea is to use the fact that a regular $n$-gon, $n\ge 7$, is made up of $n$ copies of a subequilateral triangle, and to keep Figure~\ref{F2a} in mind. When $N$ tends to infinity, the polygon $\cP_N$ tends to the disk in the Hausdorff distance. According to \cite[Remark~2, p.~206]{LeWe1986}, it follows that, for all $j \ge 1$, the Neumann eigenvalue $\nu_j(\cP_N)$ tends to the Neumann eigenvalue $\nu_j(\cD)$ of the unit disk. The Neumann eigenvalues of the unit disk satisfy
\begin{equation}\label{E-polyN-2}
\nu_1(\cD) < \nu_2(\cD) = \nu_3(\cD) < \nu_4(\cD) = \nu_5(\cD) < \nu_6(\cD) < \nu_7(\cD) \cdots
\end{equation}
and are given respectively by the squares of the zeros of the derivatives of Bessel functions: $0=j'_{0,1}, j'_{1,1}, j'_{2,1}$, $j'_{0,2}$, and $j'_{3,1}$. It follows that, for $N$ large enough, the eigenvalue $\nu_6(\cP_N)$ is simple.\medskip

From now on, we assume that $N$ is sufficiently large to ensure that $\nu_6(\cP_N)$ is a simple eigenvalue. Let $u_6$ be an associated eigenfunction. \medskip

Call $A_i, 1 \le i \le N$, the vertices of $\cP_N$, so that the triangles $OA_iA_{i+1}$ are subequilateral triangles with apex angle $\frac{2\pi}{N}$. Let $\cT_N$ be the triangle $O A_1 A_2$. \medskip

Call $D_j$ the $2N$ lines of symmetry of $\cP_N$. When $N=2m$ is even, the lines of symmetry are the $m$ diagonals joining opposite vertices, and the $m$ lines joining the mid-points of opposite sides. When $N=2m+1$ is odd, the lines of symmetries are the $N$ lines joining the vertex $A_i$ to the mid-point of the opposite side. Call $D_1$ the line of symmetry passing through the first vertex. Call $D_2$ the line of symmetry such that the angle $(D_1,D_2)$ is equal to $\pi/N$. Denote the corresponding mirror symmetries by $D_1$ and $D_2$ as well. The symmetry group of the regular $N$-gon is the dihedral group $\D_N$ with presentation,
\begin{equation}\label{E-DN}
\D_N = \set{D_1,D_2 ~|~ D_1^2 = D_2^2 = 1, (D_2D_1)^N=1}\,.
\end{equation}

\begin{figure}[!hbt]
\centering
\includegraphics[scale=0.5]{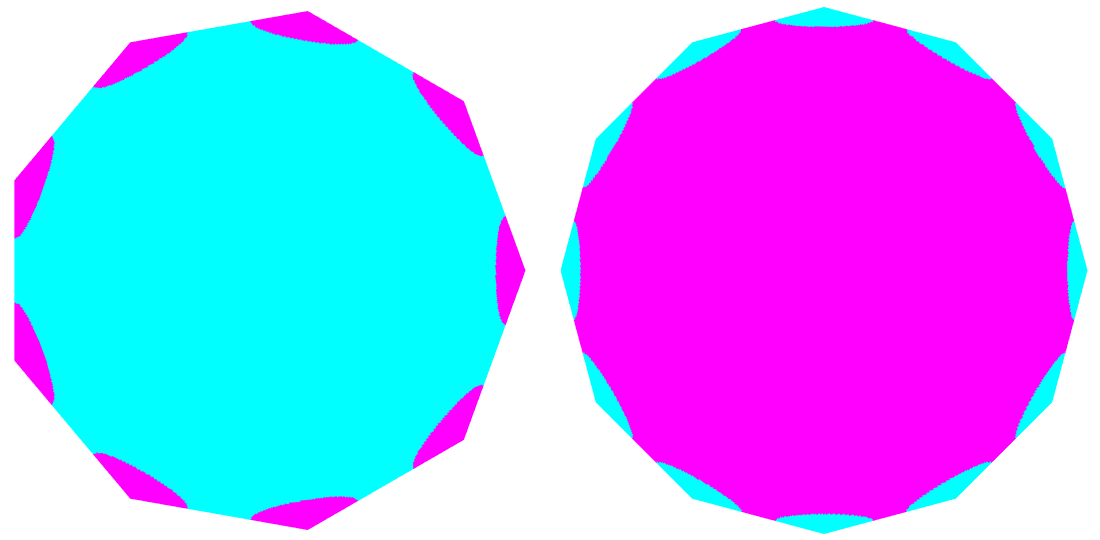}
\caption{$\cP_{9}$ and $\cP_{12}$, Neumann boundary condition}\label{F3}
\end{figure}

\FloatBarrier

The group $\D_N$ acts on functions, and commutes with the Laplacian. It leaves the eigenspaces invariant, and we therefore have a representation of degree $1$ in the eigenspace $\cE(\nu_6)$. This representation must be equivalent to one of the irreducible representations of $\D_N$ of degree $1$. When $N$ is even, there are $4$ such representations, $\rho_{\sigma_1,\sigma_2}$ with $\sigma_1,\sigma_2 \in \{-1,1\}$, and such that $\rho_{\sigma_1,\sigma_2}(D_1)=\sigma_1$ and $\rho_{\sigma_1,\sigma_2}(D_2)=\sigma_2$. When $N$ is odd, there are only $2$ irreducible representations of degree $1$, $\rho_{\sigma,\sigma}$, with $\sigma \in \{-1,1\}$. Eigenfunctions corresponding to simple eigenvalues must be invariant or anti-invariant under $D_1$ and $D_2$ depending on the signs of $\sigma_1$ and $\sigma_2$. Anti-invariant eigenfunctions must vanish on the corresponding line of symmetry. If $(\sigma_1,\sigma_2) \not = (1,1)$, the functions must have at least $N$ nodal domains. For $N \ge 7$, this is not possible for $\cE(\nu_6)$. An eigenfunction in $\cE(\nu_6)$ must be $D_1$ and $D_2$ invariant, and must therefore correspond to an eigenfunction of $\cT_N$ with Neumann boundary condition, and with eigenvalue $\nu_6 \ge \nu_2(\cT_N)$. We can now apply Proposition~\ref{P-isos-2}. This is illustrated by Figure~\ref{F3}, keeping Figure~\ref{F2a} in mind. \hfill \qed \medskip

\begin{remark}\label{R-isos-2}
The above proposition also shows that the regular $N$-gons, with the  Neumann boundary condition, provides a counterexample to Conjecture~\ref{C-ECP} and to Question~\ref{Q-ECP}(2), when $N$ is large enough.
\end{remark}%

\begin{remark}\label{R-isos-2a}
As shown in \cite{BH-ecp-bal}, Conjecture~\ref{C-ECP} is false for the regular hexagon $\cP_6$ with Neumann boundary condition. In this case, $\nu_6(\cP_6) = \nu_7(\cP_6)$, and has multiplicity $2$, with two eigenfunctions associated with different irreducible representations of $\D_6$.
\end{remark}%

\begin{remark}\label{R-isos-3}
Numerical computations indicate that the first eight Neumann eigenvalues of $\cP_7$ to $\cP_{12}$ have the same multiplicities as the first eight eigenvalues of the disk and, in particular, that $\nu_6$ is simple. Proposition~\ref{P-isos-4} is probably true for all $N \ge 6$. Numerical computations also indicate that this proposition should be true for $\cP_N$ with the Dirichlet boundary condition as well. The argument in the proof of Proposition~\ref{P-isos-4} fails in the cases $N=4$ and $N=5$ which remain open.
\end{remark}%

\section{Counterexamples on $\T^2$}\label{S-T2}

\subsection{Previous results}

This section is strongly motivated by the following result\footnote{The authors would like to thank I.~Polterovich for pointing out \cite[Section~3]{BLS}.}.

\begin{proposition}[{\cite[Section~3]{BLS}}]\label{P-BLS}
There exist a smooth metric $g$ on the torus $\mathbb T^2$, in the form  $g = Q(x) (dx^2 + dy^2)$, an infinite sequence $\phi_j$ of eigenfunctions of the Laplace-Beltrami operator $\Delta_g$, and an infinite sequence $c_j$ of real numbers, such that the level sets  $\{ (x,y) ~|~
\phi_{j}(x,y) = c_j \}$ have infinitely many connected components.
\end{proposition}

In this section, we give an easy proof of Proposition~\ref{P-BLS}, in the particular case  of one eigenfunction only, avoiding the subtleties  of \cite{BLS}. This particular case is sufficient to prove that Conjecture~\ref{C-ECP} is false on $(\T^2,g)$ for some Liouville metrics which can be chosen arbitrarily close to the flat metric.

\subsection{Metrics on $\T^2$ with a prescribed eigenfunction}\label{SS-pet2}

As in \cite{BLS}, we use the approach of Jakobson and Nadirashvili \cite{JaNa1999}. We equip the torus $\T^2 = \left( \R/ 2\pi\Z \right)^2$ with a Liouville metric of the form $g_Q = Q(x)\, g_0$, where $g_0 = dx^2 + dy^2$ is the flat metric, and where $Q$ is a positive $C^{\infty}$ function on $\T^1 = \R/ 2\pi\Z$. The respective Laplace-Beltrami operators are denoted $\Delta_0 = \partial^2_x + \partial^2_y$, and $\Delta_{Q} = Q(x)^{-1}\, \Delta_0$. \medskip

Generally speaking, we identify functions on $\T^1$ (resp. $\T^2$) with periodic functions on $\R$ (resp. $\R^2$).\medskip

Given a positive function $Q$, a complete set of spectral pairs $(\lambda,\phi)$ for the eigenvalue problem
\begin{equation}\label{E-pet2-2}
- \Delta_0 \phi(x,y) = \lambda \, Q(x)\, \phi(x,y) \text{~on~} \T^2 \,,
\end{equation}
is given by the pairs
\begin{equation}\label{E-pet2-4}
\left\{
\begin{array}{l}
\left(\sigma_{m,j},F_{m,j}(x)\cos(my)\right), m\in \N, j\in \Nb \text{~and,}\\[5pt]
\left(\sigma_{m,j},F_{m,j}(x)\sin(my)\right), m\in \Nb, j\in \Nb\,,
\end{array}%
\right.
\end{equation}
where, for a given $m \in \N$, and for $j \in \Nb$, the pair $(\sigma_{m,j},F_{m,j})$ is a spectral pair for the eigenvalue problem
\begin{equation}\label{E-pet2-1}
- u''(x) + m^2 u(x) = \sigma \, Q(x) \, u(x) \text{~on~} \T^1\,.
\end{equation}
In order to prescribe an eigenfunction, we work the other way around. Choosing a positive $C^{\infty}$ function $F$ on $\T^1$, and an integer $m \in \Nb$, we define $\Phi(x,y) = F(x) \, \cos(my)$. Then, the function $\Phi$ is an eigenfunction of the eigenvalue problem \eqref{E-pet2-2}, as soon as $Q$ and $\lambda$ satisfy,
\begin{equation*}
Q(x) = \frac{m^2}{\lambda}\, \left( 1 - \frac{1}{m^2}\, \frac{F''(x)}{F(x)}\right)\,.
\end{equation*}
Since the nodal behaviour of eigenfunctions is not affected by rescaling of the metric, we may choose $\lambda = m^2$. In order to guarantee that the metric $g_Q$ is well-defined, we need the function $Q$ to be smooth and positive. Choosing $F$ and $m$ such that
\begin{equation}\label{E-pet2-8q}
\left\{
\begin{array}{l}
0 < F(x) \text{~and,}\\[5pt]
F''(x) < m^2 \, F(x)\,, ~~\text{~for all~} x \in \R\,,
\end{array}
\right.
\end{equation}
the function $Q$, given by
\begin{equation}\label{E-pet2-8}
Q(x) = 1 - \frac{1}{m^2}\, \frac{F''(x)}{F(x)} \,,
\end{equation}
defines a Liouville metric $g_Q = Q \, g_0$ on $\T^2$, for which
$$
(\Delta_{Q}+m^2)\left( F(x)\,\cos(my)\right) =0.
$$
When $m$ is large, the metric $g_Q$ appears as a perturbation of the flat metric $g_0$.\medskip

In Subsections~\ref{SS-Ex1} and \ref{SS-Ex2}, we apply this idea to construct eigenfunctions with many level domains.

\subsection{Example~1}\label{SS-Ex1}

In this subsection, we prove the following result by describing an explicit construction.

\begin{proposition}\label{P-ex1-2}
There exists a metric $g_Q = Q(x) \left( dx^2 + dy^2\right)$ on the torus $\T^2$, and an  eigenfunction $\Phi$ of the associated Laplace-Beltrami operator, $\Delta_Q \Phi := Q^{-1}\, \Delta_0 \Phi$, such that the super-level set $\set{\Phi > 1}$ has infinitely many connected components. As a consequence, Conjecture~\ref{C-ECP} is false for $(\T^2,g_Q)$.
\end{proposition}%

\begin{remark}\label{R-ex1-2}
This proposition also implies that $\Phi$ has infinitely many isolated critical points, a particular case of \cite[Theorem~1]{BLS}.
\end{remark}%

\emph{Proof.}\\

\fbox{Step~1.}~~Let $\phi : [-\pi,\pi] \to \R$ be a function such that
\begin{equation}\label{E-ex1-2}
\left\{
\begin{array}{l}
0 \le \phi(x) \le 1\,,\\[5pt]
\supp(\phi) \subset [-\frac{\pi}{2},\frac{\pi}{2}]\,,\\[5pt]
\phi \equiv 1 \text{~on~} [-\frac{\pi}{3},\frac{\pi}{3}]\,.
\end{array}%
\right.
\end{equation}

Define the function $\psi_1 : [-\pi,\pi] \to \R$ by
\begin{equation}\label{E-ex1-4}
\psi_1(x) = \phi(x) \, \exp\left( -\frac{1}{x^2} \right)\, \cos\left( \frac{1}{x^2} \right) + 1 - \phi(x)\,.
\end{equation}

It is clear that $\psi_1$ satisfies
\begin{equation}\label{E-ex1-6}
\left\{
\begin{array}{l}
|\psi_1(x)| \le 1\,,\\[5pt]
|x| > \frac{\pi}{2} \Rightarrow \psi_1(x) = 1\,,\\[5pt]
|x| > \frac{\pi}{3} \Rightarrow \psi_1(x) > 0\,.
\end{array}%
\right.
\end{equation}

It follows that $\psi_1$ can only vanish in $[-\frac{\pi}{3},\frac{\pi}{3}]$, with zero set $\cZ(\psi_1) = \set{x ~|~ \psi_1(x)=0}$ given by
\begin{equation}\label{E-ex1-8}
\cZ(\psi_1) = \set{0} \cup \set{\pm\, \frac{1}{\sqrt{\frac{\pi}{2}+k\pi}} ~\Big| ~ k\in \N}\,.
\end{equation}

\begin{figure}[!hbt]
\centering
\includegraphics[scale=0.35]{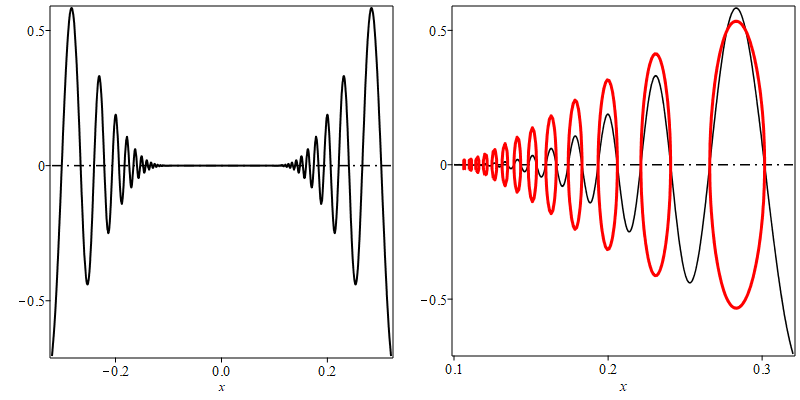}\vspace{-6mm}
\caption{Graph of $\psi_1$ near $0$. Level set $\set{\Phi = 1}$}\label{F-ex1-2}
\end{figure}

The zero set $\cZ(\psi_1)$ is an infinite sequence with $0$ as only accumulation point, and the function $\psi_1$ changes sign at each zero. The graph of $\psi_1$ over $[-\frac{\pi}{3},\frac{\pi}{3}]$ looks like the graph in the left part of Figure~\ref{F-ex1-2}.\medskip

\fbox{Step~2.}~~Define $\psi_0$ to be the function $\psi_1$ extended as a $2\pi$-periodic function on $\R$, and $F$ to be $F := 1 + \frac{1}{2} \psi_0$. Given $m\in \Nb$, define the function $\Phi_m : \T^2 \to \R$ to be $\Phi_m(x,y) = F(x) \, \cos(my)$. \medskip

The functions $F$ and $\Phi_m$ satisfy,
\begin{equation}\label{E-ex1-10}
\left\{
\begin{array}{l}
F \in C^{\infty}(\T^1)\,,\\[5pt]
F \ge \frac{1}{2}\,,\\[5pt]
\set{\psi_0 < 0}\times \T^1 \subset \set{\Phi_m < 1}\,,\\[5pt]
 \set{\Phi_m \ge 1} \subset \set{\psi_0 \ge 0} \times \T^1\,,\\[5pt]
\set{\psi_0 \ge 0}\times \set{0} \subset  \set{\Phi_m \ge 1}.
\end{array}%
\right.
\end{equation}

It follows from \eqref{E-ex1-8} that $\set{\psi_0 \ge 0} \subset \T^1$ is the union of infinitely many pairwise disjoint closed intervals, $I_{\ell}, \ell \in \Z$. It follows from \eqref{E-ex1-10} that there is at least one connected component  of the super-level set $\set{\Phi_m > 1}$ in each $I_{\ell} \times \T^1$. This construction is illustrated in Figure~\ref{F-ex1-2}. The figure on the right displays the components of $\{\Phi_1 = 1\}$ (red curves), and the part of the graph of $\psi_1$ (black curve)\footnote{As a matter of fact, we have used differently scaled functions in order to enhance the figures.} contained in $[0.1 , 0.3]\times \T^1$.  The number of connected components of $\set{\Phi_1 = 1}$ contained in $[\alpha , 0.3]\times \T^1$ tends to infinity as $\alpha$ tends to zero from above, and the components accumulate to $(0,0)$.\medskip

We have constructed a family of functions, $\Phi_m, m \in \Nb$, whose super-level sets $\set{\Phi_m > 1}$ have infinitely many connected components in $\T^2$. \medskip

\fbox{Step~3.}~~Since $F \in C^{\infty}(\T^1)$, and $F \ge \frac{1}{2}$, the function $\frac{F''}{F}$ is bounded from above. We choose $m$ such that
\begin{equation}\label{E-ex1-20}
m^2 > \sup_{x\in \T^1}\frac{F''(x)}{F(x)}\,,
\end{equation}
and we define the function $Q_m : \T^1 \to \R$,
\begin{equation}\label{E-ex1-22}
Q_m(x) = 1 - \frac{1}{m^2}\, \frac{F''(x)}{F(x)}\,.
\end{equation}

According to Subsection~\ref{SS-pet2}, and under condition \eqref{E-ex1-20}, the function $Q_m$ defines a Liouville metric $g_m$ on $\T^2$,
\begin{equation}\label{E-ex1-24}
g_m = Q_m(x)\, \left( dx^2 + dy^2 \right)
\end{equation}
and this metric can be chosen arbitrarily close to the flat metric $dx^2 + dy^2$ as $m$ goes to infinity. For the associated Laplace-Beltrami operator $\Delta_{g_m}$, we have
\begin{equation}\label{E-ex1-25}
- \Delta_{g_m} \Phi_m(x,y) = m^2 \Phi_m(x,y)\,,
\end{equation}
so that the function $\Phi_m$ is an eigenfunction of $\Delta_{g_m}$, with eigenvalue $m^2$. The super-level set $\set{\Phi_m > 1}$ has infinitely many connected components in $\T^2$. In particular, the function $\Phi_m - 1$ has infinitely many nodal domains. \hfill \qed

\begin{hide}
\begin{remark}\label{R-T2a}
One could slightly modify the above construction as follows. Replace the function $\psi_1$ by the function $\psi_{1,a} = a\, \psi_1$, where $a$ is a (small) real parameter, and extend it to a periodic function $\psi_{0,a}$. The corresponding function $F$ becomes $F_a = 1 + \frac{a}{2}\psi_{0,a}$, and the  function
\begin{equation*}
Q_{m,a}(x) = 1 - \frac{a}{m^2}\, \frac{\psi_{0,a}''(x)}{2+a\, \psi_{0,a}(x)}
\end{equation*}
defines a Liouville metric $g_{m,a}$ on $\T^2$, provided that $a$ is small enough. The function
\begin{equation*}
\Phi_{m,a}(x,y) = \left( 1 + \frac{a}{2}\, \psi_{0,a}(x) \right)\, \cos(my)
\end{equation*}
is an eigenfunction of $-\Delta_{g_{m,a}}$, associated with the eigenvalue $m^2$. The super level set $\set{\Phi_{m,a} > 1}$ has infinitely many connected components. When $m$ is fixed, and $a$ tends to zero, the metric $g_{m,a}$ tends to the metric $g_0$, and the labelling of $m^2$, as eigenvalue of $-\Delta_{g_{m,a}}$, remains bounded.
\end{remark}%
\end{hide}%

\subsection{Example~2}\label{SS-Ex2}

The metric constructed in Proposition~\ref{P-ex1-2} is smooth, not real analytic. In this subsection, we prove the following result in which we have a real analytic metric.

\begin{proposition}\label{P-ex2-2}
Let $n$ be any given integer. Then, there exists a real analytic Liouville metric $g = Q(x) \left( dx^2 + dy^2 \right)$ on $\T^2$, and an eigenfunction $\Phi$ of the associated Laplace-Beltrami operator, $-\Delta_g \Phi = \Phi$, with eigenvalue $1$, such that the super-level set $\set{\Phi > 1}$ has at least $n$ connected components. One can choose the metric $g$ arbitrarily close to the flat metric $g_0$. Taking $n \ge 4$, and $g$ close enough to $g_0$, the eigenvalue $1$ is either the second, third or fourth eigenvalue of $\Delta_g$.
\end{proposition}%

\begin{remarks}
(i)~It follows from the proposition that the function $\Phi - 1$ provides a counterexample to Conjecture~\ref{C-ECP} for $(\T^2,g)$.\\
 (ii)~The proposition is related to a question raised in \cite[Introduction]{BLS}: ``For an analytic metric, does there exist an asymptotic upper bound for the number of critical points in terms of the corresponding eigenvalue''. Indeed, given any $n \ge 4$, the function $\Phi$ given by the proposition is associated with the eigenvalue $1$, whose labelling is at most $4$, and $\Phi$ has at least $n$ isolated critical points.
 \end{remarks}

\emph{Proof.}~ Fix the integer $n$. For $0 \le a < 1$, define the functions
\begin{equation}\label{E-ex2-2}
\left\{
\begin{array}{l}
F_a(x) = 1 + a\, \cos(nx)\,,\\[5pt]
\Phi_a(x,y) = F_a(x) \, \cos(y)\,.
\end{array}%
\right.
\end{equation}

For $a$ small enough (depending on $n$), the function
\begin{equation}\label{E-ex2-4}
Q_a(x) = 1 - \frac{F_a''(x)}{F_a(x)} = 1 + a\, n^2\frac{\cos(nx)}{1 + a\, \cos(nx)}\,,
\end{equation}
is positive. According to Subsection~\ref{SS-pet2} (choosing $\psi_0(x) = a\, \cos(nx)$), the function $Q_a$ defines a Liouville metric $g_a$ on $\T^2$. The associated Laplace-Beltrami operator is $\Delta_a = \left( Q_a(x) \right)^{-1}\, \Delta_0$, and we have,
\begin{equation}\label{E-ex2-5}
\left\{
\begin{array}{l}
-\Delta_a F_a(x) \cos(y) = F_a(x) \cos(y)\,,\\[5pt]
-\Delta_a F_a(x) \sin(y) = F_a(x) \sin(y)\,.
\end{array}
\right.
\end{equation}
Call $\{\lambda_{a,j}, j\ge 1\}$ the eigenvalues of $-\Delta_a$,  written in non-decreasing order, with multiplicities.\medskip

The eigenvalues of $-\Delta_0$ are given by
\begin{equation}\label{E-ex2-6}
\left\{
\begin{array}{lll}
\lambda_{0,1} &=& 0\,,\\[5pt]
\lambda_{0,j} &=& 1  \text{~~for~~} 2 \le j \le 5 \,,\\[5pt]
\lambda_{0,j} &\ge& 2 \text{~~for~~} j\ge 6.
\end{array}
\right.
\end{equation}

For $n$ fixed, and $a$ small enough (depending on $n$), the eigenvalues $\lambda_{a,j}$ satisfy
\begin{equation}\label{E-ex2-6a}
\left\{
\begin{array}{lll}
\lambda_{a,1} &=&0\,,\\[5pt]
\lambda_{a,j} &\in& ]0.8,1.2[ \text{~~for~~} 2 \le j \le 5\,,\\[5pt]
\lambda_{a,j} &\ge& 1.8 \text{~~for~~} j\ge 6.
\end{array}
\right.
\end{equation}
We note that the metric $g_a$ and the operators $\Delta_a$ are invariant under the symmetries $\Sigma_1 : (x,y) \to (-x,y)$ and $\Sigma_2 : (x,y)\to (x, -y)$, which commute. Consequently, the space $L^2(\T^2,g_a)$ decomposes into four orthogonal subspaces
\begin{equation}\label{E-decomp}
\cS_{\varepsilon_1,\varepsilon_2} = \set{\phi \in L^2(\T^2) ~|~ \Sigma_1^*\phi = \varepsilon_1 \phi, ~\Sigma_2^*\phi = \varepsilon_2 \phi}\,,
\end{equation}
and the eigenvalue problem for $\Delta_a$ on $L^2(\T^2)$ splits into four independent problems by restriction to the subspaces $\cS_{\varepsilon_1,\varepsilon_2}$, with $\varepsilon_1,\varepsilon_2 \in \{-,+\}$. The eigenvalue $0$ is the first eigenvalue of $-\Delta_a|\cS_{+,+}$. \smallskip

When $a=0$, the eigenvalue $1$ arises with multiplicity $2$ from  $-\Delta_a|\cS_{+,+}$ (the functions $\cos x$ and $\cos y$), with multiplicity $1$ from $-\Delta_a|\cS_{-,+}$ (the function $\sin x$), and multiplicity $1$ from $-\Delta_a|\cS_{+,-}$ (the function $\sin y$).\smallskip

For $a$ small enough, the same spaces yield the eigenvalues $\lambda_{a,j}, 2 \le j \le 5$. According to \eqref{E-ex2-5}, the functions $F_a(x) \cos(y) \in \cS_{+,+}$ and $F_a(x) \sin(y) \in \cS_{+,-}$ correspond to the eigenvalue $1$. In view of \eqref{E-ex2-6a}, there is another eigenvalues $\sigma(a)$ of  $-\Delta_a|\cS_{+,+}$, and another eigenvalue $\tau(a)$ of  $-\Delta_a|\cS_{-,+}$, with $\sigma(a), \tau(a) \in \, ]0.8,1.2[$ (these eigenvalues could possibly be equal to $1$). It follows that
\begin{equation*}
\set{\lambda_{a,j}, 2 \le j \le 5} = \set{1,1,\sigma(a),\tau(a)}\,,
\end{equation*}
so that the eigenvalue $1$ of $-\Delta_a$ is either $\lambda_{a,2}$, $\lambda_{a,3}$, or $\lambda_{a,4}$. Note that, in view of \eqref{E-ex2-6a}, these eigenvalues are the smallest nonzero eigenvalues of $-\Delta_a$ restricted to the corresponding symmetry spaces $\cS$.\smallskip

Letting $\psi_0(x) = a\, \cos(nx)$, and choosing $a$ small enough (depending on $n$), the functions $F_a$ and $\Phi_a$ satisfy relations similar to the relations \eqref{E-ex1-10} of Subsection~\ref{SS-Ex1},
\begin{equation}\label{E-ex2-10}
\left\{
\begin{array}{l}
F_a \in C^{\infty}(\T^1)\,,\\[5pt]
F_a \ge \frac{1}{2}\,,\\[5pt]
\set{\psi_0 < 0}\times \T^1 \subset \set{\Phi_a < 1}\,,\\[5pt]
 \set{\Phi_a \ge 1} \subset \set{\psi_0 \ge 0} \times \T^1\,,\\[5pt]
\set{\psi_0 \ge 0}\times \set{0} \subset  \set{\Phi_a \ge 1}.
\end{array}%
\right.
\end{equation}

\begin{figure}[!hbt]
\centering
\includegraphics[scale=0.4]{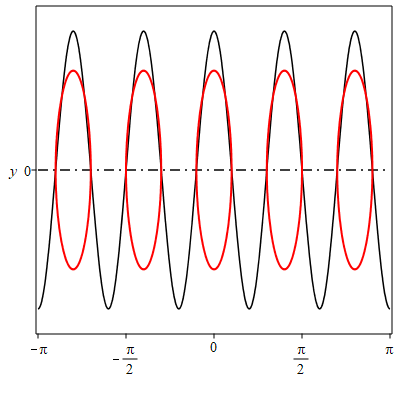}\vspace{-5mm}
\caption{Level sets $\set{\Phi_a = 1}$ for $n=5$}\label{F-ex2-2}
\end{figure}

These relations show that the super-level set $\set{\Phi_a > 1}$ has $n$ connected components. This is illustrated in Figure~\ref{F-ex2-2} for $n=5$. The black curve is the graph of $0.6\, \cos(nx)$. The red curves are components of the corresponding level set $\set{\Phi_a = 1}$  (in the figure, $a=0.1$).

It follows that the function $\Phi_a - 1$ has at least $n+1$ nodal domains. When $n \ge 4$, this also tells us that $\Phi_a-1$ provides a counterexample to Conjecture~\ref{C-ECP}.  \hfill \qed \medskip

\FloatBarrier

\subsection{Perturbation theory}\label{SS-perturb}

 We use the same notation as in Subsection~\ref{SS-Ex2}. Using perturbation theory, we now analyze the location of the eigenvalue $1$ in the spectrum of the operator $\Delta_a$, and refine Proposition~\ref{P-ex2-2}. More precisely, we prove

\begin{proposition}\label{P-perturb-2} For any given $n\geq 3$, and $a$ small enough (depending on $n$), the eigenvalue $1$ is the fourth eigenvalue of the operator $-\Delta_a$ associated with the metric $g_a = Q_a(x) \left( dx^2+dy^2\right)$, where $Q_a(x)$ is defined in \eqref{E-ex2-4}
\end{proposition}

\proof We have constructed the metrics $g_a$ in such a way that $1$ is always an eigenvalue of multiplicity at least $2$ (see \eqref{E-ex2-5}). We may assume that $a$ is small enough so that \eqref{E-ex2-6a} holds. \medskip

The idea of the proof is to show that the eigenvalues $\sigma(a)$ and $\tau(a)$ are less than $1$ by looking at their expansions\footnote{For the existence of such expansions, we refer to \cite{Rel1969} or \cite{Kat1980}.} in powers of $a$. It will actually be sufficient to compute the first three terms of these expansions. As in the proof of Proposition~\ref{P-ex2-2}, we use the symmetry properties of the metrics $g_a$.\medskip

From the proof of Proposition~\ref{P-ex2-2}, we know that $\sigma(a)$ is an eigenvalue of $-\Delta_a|\cS_{+,+}$ and $\tau(a)$ an eigenvalue of  $-\Delta_a|\cS_{-,+}$. It therefore suffices to look at eigenfunctions which are even in the variable $y$. Using Fourier cosine decomposition in the variable $y$, we reduce our problem to analyzing the family of eigenvalue problems,
\begin{equation}\label{E-pert-k}
- u''(x) + k^2 u(x) = \sigma \, Q_a(x) \, u(x) \text{~on~} \T^1\,, \text{~~~for~} k \in \N\,.
\end{equation}

More precisely, to study $\sigma(a)$, we look at the family \eqref{E-pert-k} restricted to \emph{even} functions in the $x$-variable; to study $\tau(a)$, we look at the family \eqref{E-pert-k} restricted to \emph{odd} functions in the $x$-variable.\medskip

In the sequel, we use the notation $\ip{f}{g} = \int_{-\pi}^{\pi}f(t) g(t)\, dt$ to denote the inner product of real functions in $L^2(\T^1)$.\medskip

{\bf Claim~1.}~ For $n \ge 3$, and $a$ small enough, $\sigma(a) < 1$.\medskip

Recall that to analyse $\sigma(a)$, we restrict \eqref{E-pert-k} to \emph{even functions} of $x$.\smallskip

The eigenvalue $0$ of $-\Delta_a$ appears as the first eigenvalue of \eqref{E-pert-k}, for $k=0$. When $a=0$, the eigenvalue $1$ appears as the second eigenvalue of \eqref{E-pert-k}, for $k=0$, and as the first eigenvalue of \eqref{E-pert-k}, for $k=1$. When $a>0$, the function $F_a$, defined in \eqref{E-ex2-2}, is an eigenfunction of \eqref{E-pert-k}, for $k=1$, and $1$ is the first eigenvalue of this equation because $F_a$ is positive. For $a$ small enough, the second eigenvalue of this equation must be larger than $2$. It follows from \eqref{E-ex2-6a} that $\sigma(a)$ cannot be an eigenvalue of \eqref{E-pert-k}, for $k=1$. By the min-max, $\sigma(a)$ cannot either be an eigenvalue of \eqref{E-pert-k}, for $k\ge 2$. It again follows from \eqref{E-ex2-6a} that $\sigma(a)$ must be the second eigenvalue of \eqref{E-pert-k}, for $k=0$.\medskip

We rewrite \eqref{E-pert-k}, for $k=0$, \emph{restricted to even functions}, as
\begin{equation}\label{E-pet2-1a}
- u''(x) = \sigma \, Q_a(x)  \, u(x) \text{~on~} \T^1\,, ~~~u \text{~~even.}
\end{equation}
Since $\sigma(a)$ is a simple eigenvalue of \eqref{E-pet2-1a}, the perturbative analysis in $a$ is easy. There exist expansions of the eigenvalue $\sigma(a)$ of \eqref{E-pet2-1a}, and of a corresponding eigenfunction $u(\cdot,a)$, in the form,
$$
\sigma(a)  = \sigma_0 + \sum_{j>0} \sigma_j \, a^j   \,,\, u(\cdot,a) =  u_0 + \sum_{j>0} u_j \, a^j \,,
$$
with $\sigma_0 = 1$, $u_0(x) = \cos(x)$ (respectively the unperturbed eigenvalue and eigenfunction), and with the additional orthogonality condition,
\begin{equation}\label{E-oce-2}
\ip{u(\cdot,a)}{u(\cdot,a)} = \ip{u_0}{u_0} ~~\text{for all~~} a\,.
\end{equation}

In order to prove Claim~1, it suffices to show that $\sigma_1=0$ and $\sigma_2 < 0$. For this purpose, we now determine $\sigma_1$, $\sigma_2$, and $u_1$. Developing the left-hand side of \eqref{E-oce-2} with respect to $a$, we find a series of orthogonality conditions on the functions $u_j$. In order to determine $u_1$, we only need the orthogonality condition,
\begin{equation}\label{E-oce-4}
\ip{u_0}{u_1} = 0\,.
\end{equation}

Develop the function $Q_a$ in powers of $a$,
\begin{equation}\label{E-Qexp}
Q_a(x) = 1 +n^2  \sum_{j>0}   (-1)^{j-1} \cos^j(nx) \, a ^j   \,.
\end{equation}%

Plugging the expansions of $\sigma(a), u(\cdot,a)$ and $Q_a$ into equation \eqref{E-pet2-1a}, and equating the terms in $a^k, k\ge 0$, we find equations satisfied by the functions $u_k, k\ge 0$, in the form
\begin{equation}\label{E-u-k}
- u_k'' = u_k + f_k
\end{equation}
where $f_0=0$, and for $k\ge 1$, $f_k$ depends on the $\sigma_j$, for $j\le k$, and on the $u_j$, for $j\le (k-1)$. Recall that the functions $u_k, k\ge 1$ are even and that they satisfy the orthogonality relations given by \eqref{E-oce-2}. In order to prove Claim~1, we only need to write equation \eqref{E-u-k} for $u_1$ and $u_2$, together with the associated parity and orthogonality conditions,
\begin{equation}\label{E-u-1}
\left\{
\begin{array}{l}
- u_1''(x) = u_1(x) + \sigma_1\, \cos(x) + n^2 \cos(x) \, \cos(nx)\,,\\[5pt]
\hspace{0mm} u_1 ~~\text{even,}\\[5pt]
\hspace{0mm} 0=\ip{u_1}{u_0} = \int_{-\pi}^{\pi}u_1(t) \cos(t)\, dt\,.
\end{array}
\right.
\end{equation}
and
\begin{equation}\label{E-u-2}
\left\{
\begin{array}{l}
- u_2''(x) = u_2(x)  + \sigma_2\, \cos(x) - n^2 \cos(x) \, \cos^2(nx)  \\[5pt]
\hspace{10mm} + n^2 \cos(nx) \, (u_1(x) + \sigma_1 \, \cos(x))+ + \sigma_1 \, u_1(x)\,,\\[5pt]
\hspace{0mm} u_2 ~~\text{even,}\\[5pt]
\hspace{0mm} 0=2 \ip{u_2}{u_0} + \ip{u_1}{u_1}\,.
\end{array}
\right.
\end{equation}

Taking the $L^2$ inner product of the differential equation in \eqref{E-u-1}  with $\cos x$, and for $n\ge 3$, we obtain that
 \begin{equation}
 \sigma_1=0 \,.
 \end{equation}
Since $\sigma_1=0$, the differential equation satisfied by $u_1$ in \eqref{E-u-1} becomes,
\begin{equation}\label{E-u-1a}
- u_1''(x) = u_1(x)  + n^2  \cos (nx) \cos x  \,,
\end{equation}
with $u_1$ even, satisfying \eqref{E-oce-4}. Writing
\begin{equation*}
\begin{array}{ll}
\cos (nx) \cos x  = \frac{1}{2} \left\{ \cos((n+1)x) + \cos((n-1)x) \right\},
\end{array}%
\end{equation*}
it is easy to check that the function $p(x)$ defined by
\begin{equation}\label{E-u-1b}
p(x) = \frac{n}{2}\, \left( \frac{\cos((n+1)x)}{n+2} + \frac{\cos((n-1)x)}{n-2}\right)\,,
\end{equation}
is a particular solution of this differential equation. The general solution is given by $\alpha \cos(x) + \beta \sin(x) + p(x)$, and since $u_1$ is even and orthogonal to $\cos(x)$, we find that $u_1 = p$.\medskip

Taking the fact that $\sigma_1 =0$ into account, the differential equation for $u_2$ in \eqref{E-u-2} becomes
\begin{equation}\label{E-u-2a}
-u''_2(x) =u_2(x) + n^2 \cos (nx) u_1(x) - n^2 \cos^2 (nx) \cos x + \sigma_2  \cos x \,.
\end{equation}
Taking the scalar product with $\cos x$, we obtain
$$
n^2 \int \cos^2 (nx) (\cos x)^2 dx -  n^2  \int \cos (nx) u_1(x) \cos x\,dx = \sigma_2  \int (\cos x) ^2 dx\,.
$$
The sign of  $\sigma_2$ is the same as the sign of
$$A_n:=  \int_{-\pi}^{+\pi}  \cos (nx)^2 (\cos x)^2 dx -   \int_{-\pi}^{+\pi} \cos (nx) u_1(x) \cos x\,dx \,.
$$
Since $u_1 = p$, computing each term of the sum, we get
$$
\int \cos (nx)^2 (\cos x)^2 dx
= \frac{\pi}{2}\,,
$$
and
$$
\int \cos (nx) u_1(x) \cos x\,dx
= \frac{\pi}{2} \, \frac{n^2}{n^2-4}
$$
Finally,
$$
\frac 2 \pi A_n = - \, \frac{4}{n^2-4} < 0\,.
$$
Claim~1 is proved: $\sigma(a) < 1$ for $a$ small enough.

\begin{remark}
We could continue the construction at any order, but we do not need it for our purposes.
\end{remark}%

{\bf Claim~2.}~ For $n \ge 3$ and $a$ small enough, $\tau(a) < 1$. \medskip

Recall, from the proof of Proposition~\ref{P-ex2-2}, that $\tau(a)$ is the first eigenvalue of $-\Delta_a|\cS_{-,+}$. In order to analyse $\tau(a)$, we restrict \eqref{E-pert-k} to \emph{odd functions} of $x$.\smallskip

The eigenvalue $\tau(a)$ is actually the first eigenvalue of \eqref{E-pert-k}, for $k=0$, restricted to odd functions of $x$. We rewrite \eqref{E-pert-k}, for $k=0$, \emph{restricted to odd functions}, as
\begin{equation}\label{E-pet2-1ao}
- u''(x) = \sigma \, Q_a(x)  \, u(x) \text{~on~} \T^1\,, ~~~u \text{~~odd.}
\end{equation}
Since $\tau(a)$ is a simple eigenvalue of \eqref{E-pet2-1ao}, the perturbative analysis in $a$ is easy. There exist expansions of the eigenvalue $\tau(a)$ of \eqref{E-pet2-1ao}, and of a corresponding eigenfunction $v(\cdot,a)$, in the form,
$$
\tau(a)  = \tau_0 + \sum_{j>0} \tau_j \, a^j   \,,\, v(\cdot,a) =  v_0 + \sum_{j>0} v_j \, a^j \,,
$$
with $\tau_0 = 1$, $v_0(x) = \sin(x)$ (respectively the unperturbed eigenvalue and eigenfunction), and with the additional orthogonality condition,
\begin{equation}\label{E-oce-2o}
\ip{v(\cdot,a)}{v(\cdot,a)} = \ip{v_0}{v_0} ~~\text{for all~~} a\,.
\end{equation}

To prove Claim~2, it suffices to show that  $\tau_1=0$ and $\tau_2 < 0$. For this purpose, it is sufficient to determine $\tau_1$, $\tau_2$, and $v_1$. Developing the left-hand side of \eqref{E-oce-2o} with respect to $a$, we find a series of orthogonality conditions on the functions $v_j$. In order to determine $v_1$, we only need the orthogonality condition,
\begin{equation}\label{E-oce-4o}
\ip{v_0}{v_1} = 0\,.
\end{equation}

Plugging the expansions of $\tau(a), v(\cdot,a)$ and $Q_a$, see \eqref{E-Qexp}, into equation \eqref{E-pet2-1ao}, and equating the terms in $a^k, k\ge 0$, we find equations satisfied by the functions $v_k, k\ge 0$, in the form
\begin{equation}\label{E-u-ko}
- v_k'' = v_k + h_k
\end{equation}
where $h_0=0$, and for $k\ge 1$, $h_k$ depends on the $\tau_j$, for $j\le k$, and on the $v_j$, for $j\le (k-1)$. Recall that the functions $v_k, k\ge 1$ are odd and that they satisfy the orthogonality relations given by \eqref{E-oce-2o}. In order to prove Claim~2, we only need to write equation \eqref{E-u-ko} for $v_1$ and $v_2$, together with the associated parity and orthogonality conditions,
\begin{equation}\label{E-u-1o}
\left\{
\begin{array}{l}
- v_1''(x) = v_1(x) + \tau_1\, \sin(x) + n^2 \sin(x) \, \cos(nx)\,,\\[5pt]
\hspace{0mm} v_1 ~~\text{odd,}\\[5pt]
\hspace{0mm} 0=\ip{v_1}{v_0} = \int_{-\pi}^{\pi}v_1(t) \sin(t)\, dt\,.
\end{array}
\right.
\end{equation}
and
\begin{equation}\label{E-u-2o}
\left\{
\begin{array}{l}
- v_2''(x) = v_2(x)  + \tau_2\, \sin(x) - n^2 \sin(x) \, \cos^2(nx)  \\[5pt]
\hspace{10mm} + n^2 \cos(nx) \, (v_1(x) + \tau_1 \, \sin(x)) + \tau_1 \, v_1(x)\,,\\[5pt]
\hspace{0mm} v_2 ~~\text{odd,}\\[5pt]
\hspace{0mm} 0=2 \ip{v_2}{v_0} + \ip{v_1}{v_1}\,.
\end{array}
\right.
\end{equation}

Taking the $L^2$ inner product of the differential equation in \eqref{E-u-1o}  with $\sin x$, and assuming that $n\ge 3$, we obtain that
 \begin{equation}
 \tau_1=0 \,.
 \end{equation}
Since $\tau_1=0$, the differential equation satisfied by $v_1$ in \eqref{E-u-1o} becomes,
\begin{equation}\label{E-u-1ao}
- v_1''(x) = v_1(x)  + n^2  \cos (nx) \sin(x)  \,,
\end{equation}
with $v_1$ odd, satisfying \eqref{E-oce-4o}. Writing
\begin{equation*}
\begin{array}{ll}
\cos (nx) \sin x  = \frac{1}{2} \left\{ \sin((n+1)x) - \sin((n-1)x) \right\},
\end{array}%
\end{equation*}
it is easy to check that the function $q(x)$ defined by
\begin{equation}\label{E-u-1bo}
q(x) = \frac{n}{2}\, \left( \frac{\sin((n+1)x)}{n+2} - \frac{\sin((n-1)x)}{n-2}\right)\,,
\end{equation}
is a particular solution of this differential equation. The general solution is given by $\alpha \cos(x) + \beta \sin(x) + q(x)$, and since $v_1$ is odd and orthogonal to $\sin(x)$, we find that $v_1 = q$.\medskip

Taking the fact that $\tau_1 =0$ into account, the differential equation for $v_2$ in \eqref{E-u-2o} becomes
\begin{equation}\label{E-u-2ao}
-v''_2(x) =v_2(x) + n^2 \cos (nx) v_1(x) - n^2 \cos^2 (nx) \sin x + \tau_2  \sin x \,.
\end{equation}
Taking the scalar product with $\sin x$, we obtain
$$
n^2 \int \cos^2 (nx) \sin^2( x)\, dx -  n^2  \int \cos (nx) v_1(x) \sin x\,dx = \tau_2  \int \sin^2(x)\, dx\,.
$$
The sign of  $\tau_2$ is the same as the sign of
$$B_n:=  \int_{-\pi}^{+\pi}  \cos (nx)^2 (\sin x)^2 dx -   \int_{-\pi}^{+\pi} \cos (nx) v_1(x) \sin x\,dx \,.
$$
Computing each term of the sum, we get
$$
\int \cos (nx)^2 (\sin x)^2 dx
= \frac{\pi}{2}\,,
$$
and
$$
\int \cos (nx) v_1(x) \sin x\,dx
= \frac{\pi}{2} \, \frac{n^2}{n^2-4}
$$
Finally,
$$
\frac 2 \pi B_n = - \, \frac{4}{n^2-4} < 0\,.
$$
Claim~2 is proved: $\tau(a) < 1$ for $a$ small enough.\medskip

It follows that $1$ is an eigenvalue of $-\Delta_a$, with multiplicity $2$ and least labelling $4$. This proves Proposition~\ref{P-perturb-2}.\hfill \qed

\subsection{Comparison with a result of Gladwell and Zhu}\label{SS-GZ}~\\
The authors of \cite{GlZh2003} prove the following result for a bounded domain in $\R^d$.

\begin{proposition}\label{P-GZ-2}
Let $\Omega \subset \R^d$ be a connected bounded domain. Call $(\delta_j,u_j)$ the eigenpairs of the Dirichlet eigenvalue problem in $\Omega$,
\begin{equation}\label{E-GZ-2}
\left\{
\begin{array}{rr}
-\Delta u &= \delta\,  u \text{~in~} \Omega\,,\\[5pt]
u &= 0 \text{~on~} \partial \Omega\,,
\end{array}%
\right.
\end{equation}
where the eigenvalues $\delta_1 < \delta_2 \le \delta_3 \le \dots$ are listed in non-decreasing order, with multiplicities. Assume that the first eigenfunction $u_1$  is positive. For $n \ge 2$, let $v = u_n + c u_1$, for some positive constant $c$. Then, the function $v$ has at most $(n-1)$ positive sign domains, i.e., the super-level set $\set{v > 0}$ has at most $(n-1)$ connected components.
\end{proposition}%

The same   result  is true if instead of the Dirichlet boundary condition, one considers the Neumann boundary condition (assuming in this case that $\partial \Omega$ is smooth enough), or if one considers a closed real analytic Riemannian surface\footnote{It might be necessary to use a real analytic surface in order to apply Green's theorem to the nodal sets of a linear combination of eigenfunctions.}.\medskip

A more convenient formulation, is as follows. For a function $w$, and $\varepsilon \in \{-,+\}$, define $\beta_0^{\, \varepsilon}(w)$ to be the number of nodal domains of $w$, on which $\varepsilon \, w$ is positive. Proposition~\ref{P-GZ-2} can be restated as follows. For any $n \ge 2$, and any real nonzero constant $c$,
\begin{equation}\label{E-GZ-4}
\beta_0^{\, \sign(c\, u_1)}(u_n+cu_1) \le (n-1)\,.
\end{equation}

Proposition~\ref{P-GZ-2} is weaker than Conjecture~\ref{C-ECP}. Indeed, it only gives control on the number of nodal domains where the function $u_n + cu_1$ has the sign of $\sign(cu_1)$. Propositions~\ref{P-ex1-2} and \ref{P-ex2-2} show that one can a priori not control $ \beta_0^{\, -\sign(cu_1)}(w)$, at least in the case of the Neumann (or empty) boundary condition. However, one can observe that, fixing $n_0$, it is easy to construct  examples for which Conjecture~\ref{C-ECP} is true for all linear combinations of the $n$ first eigenfunctions,  $w \in \cL_{n}$, with $n \le n_0$~. Indeed, for $L$ large, the rectangle $]0,1[\times ]0,L[$ provides such an example for the Dirichlet or Neumann boundary conditions. More generally, one can consider manifolds which collapse on a lower dimensional manifold for which the Extended Courant property is true. 


\section{Counterexamples on $\Ss^2$}\label{S-S2m}

\subsection{Results and general approach}\label{SS-S2mr}

In this section, we extend, to the case of the sphere, the construction made in Section~\ref{S-T2}. We prove the following results.

\begin{proposition}\label{P-S2m-2}
There exist $C^{\infty}$ functions $\Phi$ and $G$ on $\Ss^2$, with the following properties.
\begin{enumerate}
  \item The super-level set $\set{\Phi > 1}$ has infinitely many connected components.
  \item The function $G$ is positive, and defines a conformal metric\break  $g_{G} = G\, g_0$ on $\Ss^2$ with associated Laplace-Beltrami operator \break $\Delta_{G} = G^{-1}\, \Delta_0$.
  \item $- \Delta_{G} \Phi = 2\, \Phi$.
  \item The eigenvalue $2$ of $-\Delta_G$ has labelling at most $4$.
\end{enumerate}
\end{proposition}%

\begin{proposition}\label{P-S2m-8}
There exists $M > 0 $ such that, for any $m\ge M$, there exist $C^{\infty}$ functions $\Phi_m$ and $G_m$ on $\Ss^2$ with the following properties.
\begin{enumerate}
  \item The super-level set $\set{\Phi_{m} > 1}$ has infinitely many connected components.
  \item The function $G_{m}$ is positive, and defines a conformal metric $g_{m} = G_{m}\, g_0$ on $\Ss^2$ with associated Laplace-Beltrami operator $\Delta_{G_m} = G_m^{-1}\, \Delta_0$.
  \item For $m \ge M$, $\left(1 - \frac{2}{m}\right) \le G_m \le \left(1 + \frac{2}{m}\right)$, and
  \item $- \Delta_{G_m} \Phi_m = m(m+1)\, \Phi_m$.
\end{enumerate}
\end{proposition}%

These propositions provide counterexamples to Conjecture \ref{C-ECP} and to Questions~\ref{Q-ECP} on the sphere.

\begin{remark}
The eigenfunctions on $\Ss^2$ constructed in the above propositions have infinitely many isolated critical points. For a similar result on $\T^2$, see Remark~\ref{R-ex1-2}(ii) which is a particular case of \cite[Theorem~1]{BLS}.
\end{remark}%

The approach is inspired by Section~\ref{S-T2}, with the following steps.
\begin{enumerate}
  \item Start from a special spherical harmonic $Y$ of the standard sphere $(\Ss^2,g_0)$, with eigenvalue $m(m+1)$.
  \item Modify $Y$ into a smooth function $\Phi$, whose super-level set $\set{\Phi > 1}$ has infinitely many connected components.
  \item Construct a conformal metric $g_Q = Q \, g_0$ on $\Ss^2$, whose associated Laplace-Beltrami operator has $\Phi$ as eigenfunction, with eigenvalue $m(m+1)$.
\end{enumerate}

The proof of Propositions~\ref{P-S2m-2} and \ref{P-S2m-8}, following the above steps, is split into the next subsections.

\subsection{Metrics on $\Ss^2$ with a prescribed eigenfunction} \label{SS-pe}

Let $g_0$ be the standard metric on the sphere
$$
\mathbb{S}^2 = \set{(x,y,z) \in \R^3 ~|~ x^2+y^2+z^2 =1}.
$$
The spherical coordinates are $(\theta, \phi) \mapsto \big( \sin\theta \, \cos\phi, \sin\theta \, \sin\phi,  \cos\theta\big)$, with $(\theta,\phi) \in ]0,\pi[\times]0,2\pi[$. In these coordinates, $$g_0 = d\theta^2 + \sin^2\theta \, d\phi^2\,,$$ the associated measure is $\sin\theta \, d\theta\, d\phi$, and the Laplace-Beltrami operator of $g_0$ is given by
$$
\Delta_0 = \frac{1}{\sin \theta}\, \frac{\partial}{\partial \theta}\big( \sin\theta \, \frac{\partial}{\partial \theta}\big) +  \frac{1}{\sin^2 \theta}\, \frac{\partial^2}{\partial \phi^2}\,.
$$

We consider conformal metrics on $\Ss^2$, in the form $g_Q = Q \, g_0$, where $Q$ is $C^{\infty}$ and positive. We denote the associated Laplace-Beltrami operator by
$$\Delta_Q =Q^{-1}\,\Delta_0\,.
$$
We assume that $Q$ is invariant under the rotations with respect to the $z$-axis, i.e., that $Q$ only depends on the variable $\theta$. \medskip

Let $\Phi$ be a smooth function on $\Ss^2$, given in spherical coordinates by $\Phi(\theta,\phi) = T(\theta)P(\phi)$. If $\Phi$ is an eigenfunction of $-\Delta_Q$
associated with the eigenvalue $\lambda$, then the functions $T$ and $P$ satisfy the equations,
\begin{align}
P''(\phi) + m^2 P(\phi)=0\,,\\
\sin(\theta)\left( \sin(\theta )T'(\theta) \right)' + \left( \lambda Q(\theta) \sin^2(\theta)-m^2\right) \, T(\theta) = 0\,,
\end{align}
where $m$ is an integer. When $Q\equiv1$, the solutions are the spherical harmonics of degree $m$, $Y_m^k, -m \le k \le m$ (as given for example in \cite[p.~302]{Ley1996}).\medskip

For $m \ge 1$, we consider the special spherical harmonic
$$
Y_{m}^m(\theta,\phi) = \sin^{m}(\theta) \,\cos(m\phi)\,.
$$
We could consider $\sin^{m}(\theta)\,\sin(m\phi)$ as well, since $m \geq 1$. For later purposes, we introduce the linear differential operator $\cK_{m}$, defined by
\begin{equation}\label{E-Km}
\begin{array}{ll}
T\mapsto ( \cK_{m}T)(\theta) &= \sin^2(\theta)\,T''(\theta) + \sin(\theta)\cos(\theta)\, T'(\theta)\\[5pt]
&~~\qquad + \left( m(m + 1) \sin^2(\theta) - m^2\right)\, T(\theta) \,.
\end{array}%
\end{equation}

In particular, we have
\begin{equation}\label{E-SHm}
\cK_m \sin^m(\cdot) = 0.
\end{equation}
Given $Q$ a smooth positive function, which only depends on $\theta$, a ne\-ces\-sary and sufficient condition for the function $\Phi(\theta,\phi) = T(\theta)\,\cos(m\phi)$ to satisfy $- \Delta_Q \Phi = m(m+1)\, \Phi$, is that
\begin{equation}\label{E-QQm}
Q(\theta) = -\, \frac{\sin^2(\theta)T''(\theta) + \sin(\theta)\cos(\theta)T'(\theta)  - m^2 T(\theta) }{m(m+1)\sin^2(\theta) T(\theta)}\,,
\end{equation}
or, equivalently,
\begin{equation}\label{E-QQmm}
1 - Q(\theta) = \frac{(\cK_mT) (\theta)}{m(m+1)\sin^2(\theta) T(\theta)}
\end{equation}

In particular, taking $\Phi(\theta,\phi) = \sin^m(\theta)\, \cos(m\phi)$, we find that $Q \equiv 1$.

\begin{remark}
As in Section~\ref{S-T2}, we work the other way around: we prescribe $T$, and look for a conformal metric on $\bS^2$ admitting $\Phi(\theta,\phi)=T(\theta)\cos(m\phi)$ as eigenfunction. The main difficulty in prescribing the function $T$, is to show that the function $Q$ defined by \eqref{E-QQm} is actually smooth and positive.
\end{remark}%

\subsection{Constructing perturbations of the function $\sin^m(\theta)$} \label{SS-psm}~\\

In Section~\ref{S-T2}, we perturbed the eigenfunction $\cos(my)$ of the torus into the function $\Phi_m(x,y) = F(x) \cos(my)$, where $F$ had rapidly decaying oscillations around $x = 0$. We do a similar construction here, with an extra flattening step. \medskip

Given $m \ge 1$, we start from the spherical harmonic $\sin^m(\theta)\, \cos(m\phi)$ in spherical coordinates. We first flatten the function $\sin^m(\theta)$ around \break $\theta = \pi/2$, before adding the rapidly decaying oscillations. More precisely, we look for functions $\Phi$ of the form $\Phi(\theta,\phi) = T(\theta) \cos(m\phi)$. To determine $T$, we construct a family $T_{m,n,\alpha}$ of perturbations of the function $\sin^m(\cdot)$, in the form,
\begin{equation}\label{E-S2m-8}
T_{m,n,\alpha}(\theta) = \sin^m(\theta) + P_{m,n,\alpha}(\theta) + u_{m,n,\alpha}(\theta)\,,
\end{equation}
with $n\in \N$ (to be chosen large), and $\alpha \in (0,\frac 14]$ (to be chosen small). The function $P_{m,n,\alpha}$ is constructed such that
\begin{equation}\label{E-S2m-P1}
\sin^m(\theta) + P_{m,n,\alpha}(\theta) \equiv 1
\end{equation}
in an interval around $\frac{\pi}{2}$, and $u_{m,n,\alpha}$ is a rapidly oscillating function in the same interval. They will both be designed in such a way that we can control the derivatives in equation (\ref{E-QQm}). The construction of the family $T_{m,n,\alpha}$ is explained in the following paragraphs, and illustrated in Figure~\ref{F-sph-2}.\medskip

\subsubsection{Construction of $P_{m,n,\alpha}$}

\begin{proposition}[Construction of $P_{m,n,\alpha}$]\label{P-S2-Pmn} For all $m \ge 1$ and $\alpha \in (0,\frac 14]$, there exist $N(m,\alpha) \in \N$, and a sequence of functions $(P_{m,n,\alpha})_{n\ge 1}$, $P_{m,n,\alpha} : [0,\pi]  \to \R$, with the following properties for all $n \ge N(m,\alpha)$.
\begin{enumerate}
	\item  $P_{m,n,\alpha} \in C^{\infty}$  and  $P_{m,n,\alpha}(\pi - \theta) = P_{m,n,\alpha}(\theta)$ for all $\theta \in [0,\pi]$;
	\item  for  $\theta \in [\frac{\pi}{2} + \frac{1}{mn}, \pi]$, $P_{m,n,\alpha}(\theta) = 0$;
	\item for  $\theta \in [\frac{\pi}{2}, \frac{\pi}{2} + \frac{\alpha}{(mn)^2}]$, $P_{m,n,\alpha}(\theta) = 1-\sin^m(\theta)$;
	\item for $\theta \in [\frac{\pi}{2}, \pi]$, $0 \le P_{m,n,\alpha}(\theta) \leq \frac{2m}{(mn)^3}$, and $|P'_{m,n,\alpha}(\theta)| \leq \frac{2m}{(mn)^2}$;
	\item for $\theta \in [\frac{\pi}{2}, \pi]$, $ -m(1 +5 \alpha) \leq P''_{m,n,\alpha}(\theta) \leq m(1 + \alpha)$.
\end{enumerate}
 \end{proposition}%

The idea is to construct $P_{m,n,\alpha}$ as
\begin{equation}\label{E-S2-Pn-2m}
P_{m,n,\alpha}(\theta) = \int_{\frac \pi2}^{\theta} R_{m,n,\alpha}(t)\, dt\,,
\end{equation}
for $\theta \in [\frac \pi2,\pi]$, and to extend it so that $P_{m,n,\alpha}(\pi - \theta) = P_{m,n,\alpha}(\theta)$. We first construct a sequence $S_{m,n,\alpha}$ (Lemma~\ref{L-S2-Smn}), and then a sequence $s_{m,n,\alpha}$, such that $R_{m,n,\alpha} = S_{m,n,\alpha} + s_{m,n,\alpha}$ (Lemma~\ref{L-S2-Rmn}).

\begin{lemma}[Construction of $S_{m,n,\alpha}$]\label{L-S2-Smn} For any $m \ge 1$, and any $\alpha \in (0,\frac 14] $, there exists a sequence of functions $(S_{m,n,\alpha})_{n\geq 1}$, $S_{m,n,\alpha} : [0,\pi]  \to \R$, with the following properties for $n\ge 2$.
\begin{enumerate}
	\item  $S_{m,n,\alpha} \in C^{\infty}$  and  $S_{m,n,\alpha}(\pi-\theta) = - S_{m,n,\alpha}(\theta)$ for all $\theta \in [0,\pi]$;
	\item  for $\theta \in [\frac{\pi}{2} + \frac{1}{(mn)^2}, \pi]$, $S_{m,n,\alpha}(\theta) = 0\,$;
	\item for  $\theta \in [\frac{\pi}{2}, \frac{\pi}{2} + \frac{\alpha}{(mn)^2}]$, $S_{m,n,\alpha}(\theta) = -m\, \cos(\theta)\,\sin^{m-1}(\theta)\,$;
	\item for $\theta \in [\frac{\pi}{2}, \pi]$, $0 \le S_{m,n,\alpha}(\theta) \leq \frac{1}{mn^2}$\,;
	\item for $\theta \in [\frac{\pi}{2}, \pi]$, $-m(1 + 4\alpha) \leq S_{m,n,\alpha}'(\theta) \leq m$\,.
\end{enumerate}
\end{lemma}

\emph{Proof of Lemma~\ref{L-S2-Smn}.}~ We construct $S_{m,n,\alpha}$ on $[\frac{\pi}{2},\pi]$, and extend it to $[0,\pi]$ so that $S_{m,n,\alpha}(\pi-\theta) = -S_{m,n,\alpha}(\theta)$.\medskip

Choose a function $\chi_{\alpha} : \R \to [0,1]$, such that $\chi_{\alpha}$ is smooth and even, $\chi_{\alpha}(t) = 1$ on $[-\alpha,\alpha]$, $\supp(\chi_{\alpha}) \subset [-1,1]$, and
 \begin{equation}\label{E-S2m-chip}
  -1-4\alpha \le - \frac{1}{1-2\alpha} \leq  \chi'_\alpha (t) \leq 0\,,\, \forall t \ge 0\,.
 \end{equation}
A natural Lipschitz candidate would be a piecewise linear function $\xi_{\alpha}$ which is equal to $1$ in $[0,\alpha]$, and to $t\mapsto 1 - (t-\alpha)/(1-\alpha)$ in $[\alpha,1]$. To get $\chi_{\alpha}$, we can regularize a function $\xi_{\beta}$, keeping the other properties at the price of a small loss in the control of the derivative in \eqref{E-S2m-chip}.\medskip

For $\theta \in [\frac{\pi}{2},\pi]$, we introduce $\hat \theta = \theta - \frac{\pi}{2}$. \medskip

We take $S_{m,n,\alpha}$ in the form
\begin{equation}\label{E-S2m-4}
\left\{
\begin{array}{ll}
S_{m,n,\alpha}(\theta) &= -m\ \chi_\alpha ((mn)^2\, (\theta - \frac{\pi}{2}))\, \cos(\theta) \, \sin^{m-1}(\theta)\\[5pt]
&= m \, \chi_\alpha ((mn)^2 \, \hat \theta)\, \sin(\hat \theta) \, \cos^{m-1}(\hat \theta)\,.
\end{array}%
\right.
 \end{equation}

Properties (1), (2) and (3) are clear. Property (4) follows from the inequality $ |\sin(\hat \theta)| \le |\hat \theta|$ and (2). To prove (5), we introduce
$$
h(\theta) := - m\, \cos(\theta)\, \sin^{m-1}(\theta)\,.
$$
For $m\ge 3$, we have
$$
h'(\theta) = m^2 \, \sin^{m-2}(\theta) \, \left( \frac{1}{m} - \sin^2(\hat\theta)\right)\,,
$$
and hence, $$m\, \sin^{m-2}(\theta) \ge h'(\theta) \ge \frac{m}{2}\sin^{m-2}(\theta)\,,$$ in the set $\set{\theta ~|~ 0\le (mn)^2\hat\theta \le 1}$, as soon as $n \ge 2\,$.\\
 We have,
$$
 S'_{m,n,\alpha}(\theta) = \chi_\alpha ((mn)^2\hat \theta ) \, h'(\theta) + (mn)^2\, \chi_\alpha'((mn)^2\hat \theta ) \, h(\theta)\,.
$$
Using the inequality $|\cos \theta| = |\sin \hat \theta| \leq |\hat \theta|$, and \eqref{E-S2m-chip} for $0 \le \hat \theta \le \frac{1}{(mn)^2}$, we obtain
 \begin{equation}
 - m(1+4\alpha) \leq S'_{m,n,\alpha}(\theta)\leq   m\,.
 \end{equation}
One can check directly that this inequality also holds for $m=1$ and $2$.  Lemma~\ref{L-S2-Smn} is proved.  \hfill \qed \medskip

\begin{lemma}[Construction of $R_{m,n,\alpha}$]\label{L-S2-Rmn}
 For all $m \ge 1$ and $\alpha \in (0,\frac 14]$, there exist $N(m,\alpha) \in \N$, and a sequence of functions $(R_{m,n,\alpha})_{n\ge 1}$, $R_{m,n,\alpha} : [0,\pi]  \to \R$ with the following properties for all $n \ge N(m,\alpha)$.
\begin{enumerate}
	\item  $R_{m,n,\alpha} \in C^{\infty}$ and $R_{m,n,\alpha}(\pi-\theta) = - R_{m,n,\alpha}(\theta)$ for all $\theta \in [0,\pi]$;
	\item  for  $\theta \in [\frac{\pi}{2} + \frac{1}{mn}, \pi]$, $R_{m,n,\alpha}(\theta) = 0\,$;
	\item for  $\theta \in [\frac{\pi}{2}, \frac{\pi}{2} +\frac{\alpha}{(mn)^2}]$, $R_{m,n,\alpha}(\theta) = -m\, \cos(\theta)\, \sin^{m-1}(\theta)\,$;
	\item for $\theta \in [\frac{\pi}{2}, \pi]$, $|R_{m,n,\alpha}(\theta)| \leq \frac{2}{mn^2}\,$;
	\item for $\theta \in [\frac{\pi}{2}, \pi]$, $ -m(1 +5 \alpha) \leq R_{m,n,\alpha}'(\theta) \leq m(1 + \alpha)\,$;
	\item   $\int\limits_{\frac{\pi}{2}}^{\frac{\pi}{2} + \frac{1}{mn}} R_{m,n,\alpha}(\theta) d\theta =0\,$;
\item $\int\limits_{\frac{\pi}{2}}^{\frac{\pi}{2} + \frac{1}{mn}} |R_{m,n,\alpha}(\theta)| d\theta \le \frac{2m}{(mn)^3}\,$.
\end{enumerate}
\end{lemma}

\emph{Proof of Lemma~\ref{L-S2-Rmn}.}~ We construct $R_{m,n,\alpha}$ on $[\frac{\pi}{2},\pi]$, and extend it to $[0,\pi]$ so that $R_{m,n,\alpha}(\pi-\theta) = -R_{m,n,\alpha}(\theta)$. Define
\begin{equation}\label{E-S2m-beta}
\beta_{m,n,\alpha} := \int_{\frac{\pi}{2}}^{\frac{\pi}{2} + \frac{1}{(mn)^2}} S_{m,n,\alpha}(t) \, dt.
\end{equation}
Using \eqref{E-S2m-4} and a change of variable, we find that $\beta_{m,n,\alpha}$ satisfies,
\begin{equation*}
\left\{
\begin{array}{l}
\beta_{m,n,\alpha} \ge m\, \int_{0}^{\frac{\alpha}{(mn)^2}} \sin(t)\, \cos^{m-1}(t) \, dt \,, \text{and}\\[5pt]
\beta_{m,n,\alpha} \le m\, \int_{0}^{\frac{1}{(mn)^2}} \sin(t)\, \cos^{m-1}(t)\, dt\,.
\end{array}%
\right.
\end{equation*}
Using the inequalities $\frac{2}{\pi} t \le \sin(t) \le t$ for $t \in [0,\frac{\pi}{2}]$, we obtain,
\begin{equation}\label{E-S2m-10}
\frac{\alpha^2}{\pi}\, \frac{m}{(mn)^4}\le \frac{\alpha^2}{\pi}\, \frac{m}{(mn)^4}\,2 \cos^{m-1}(\frac{\alpha}{(mn)^2}) \le \beta_{m,n,\alpha} \le \frac{1}{2}\, \frac{m}{(mn)^4}\,,
\end{equation}
where the first inequality holds provided that $n$ is larger than some $N_1(m)$.\medskip

Choose a $C^{\infty}$ function $\xi$, such that $0 \le \xi \le 1$, $\supp(\xi) \subset (\frac{1}{2},1)$, and $\int_{\R} \xi(t) \, dt = 1$. Note that for $n \ge 3$, $[0,\frac{1}{(mn)^2}] \cap [\frac{1}{2mn},\frac{1}{mn}] = \emptyset$\,.\smallskip

Define $s_{m,n,\alpha}$ by,
\begin{equation}
s_{m,n,\alpha}  (\theta) = - \, \gamma_{m,n,\alpha} \; \xi (mn \, \hat \theta)\,,
\end{equation}
where $\gamma_{m,n,\alpha}$ is a constant to be chosen later. Note that for $n\ge 3$, $\supp(S_{m,n,\alpha}) \cap
\supp(s_{m,n,\alpha}) = \emptyset$\,.\medskip

Defining
\begin{equation}\label{E-S2m-R}
R_{m,n,\alpha} = S_{m,n,\alpha} + s_{m,n,\alpha}\,,
\end{equation}
Assertion~(2) is satisfied. Choosing,
\begin{equation}\label{E-S2m-gam}
\gamma_{m,n,\alpha} = mn \,\beta_{m,n,\alpha}\,,
\end{equation}
Assertion~(6) is satisfied,
\begin{equation*}
\int_{\frac{\pi}{2}}^{\frac{\pi}{2} + \frac{1}{mn}} R_{n,\alpha}(t) \, dt = \int_{\frac{\pi}{2}}^{\frac{\pi}{2} + \frac{1}{mn}} \left( S_{n,\alpha}(t) + s_{n,\alpha}(t) \right)\, dt = 0\,,
\end{equation*}
and the function $\int_{\frac{\pi}{2}}^{\theta} R_{n,\alpha}(t) \, dt$ vanishes for $\theta \ge \frac{\pi}{2} + \frac{1}{mn}$.\smallskip

Using \eqref{E-S2m-10}, we get
\begin{equation}\label{E-S2-gamm}
\frac{\alpha^2}{\pi}\, \frac{m}{(mn)^3} \le \gamma_{m,n,\alpha} \le \frac{1}{2}\, \frac{m}{(mn)^3}, \text{~for~} n \ge N_1(m)\,.
\end{equation}

Properties~(1) and (3) are clear. Using the properties of $S_{m,n,\alpha}$ given by Lemma~\ref{L-S2-Smn}, inequality \eqref{E-S2-gamm}, and the fact that $\frac{1}{2} \le mn \, \hat \theta \le 1$ when $\xi' \not = 0$, we obtain
\begin{equation*}
|R_{m,n,\alpha}(\theta)| \le \frac{2}{mn^2}\,,
\end{equation*}
and
\begin{equation*}
-m \left( 1+4\alpha + \frac{1}{2(mn)^2}\,\|\xi'\|_{\infty}\right) \le R'_{m,n,\alpha}(\theta) \le m \left( 1 + \frac{1}{2(mn)^2}\,\|\xi'\|_{\infty} \right)\,.
\end{equation*}

Assertions~(4) and (5) follow by taking $n$ larger than some $N(m,\alpha)$. Assertion~(7) follows from Property~(4). Lemma~\ref{L-S2-Rmn} is proved.\hfill \qed \medskip

\emph{Proof of Proposition~\ref{P-S2-Pmn}.~} Recall that
$$
P_{m,n,\alpha} (\theta) = \int\limits_{\pi/2}^{\theta} R_{m,n,\alpha} (z)dz \; \ge 0\,,
$$
for $\theta \geq \pi/2$, and that $P_{m,n,\alpha}$ is symmetric with respect to $\frac \pi 2$. The properties of $P_{m,n,\alpha}$ follow from Lemma~\ref{L-S2-Rmn}. For $ |\theta-\frac{\pi}{2}| \leq \frac{\alpha}{(mn)^2}$, $P_{m,n,\alpha}(\theta) = 1- \sin^{m}(\theta)$. We also note that
$P_{m,n,\alpha} (\theta) =0$ for $ \theta \in  (0, \frac{\pi}{2} -  \frac{1}{mn}) \cup (\frac{\pi}{2} + \frac{1}{mn}, \pi)$. \hfill \qed \bigskip

\begin{figure}[!hbt]
\centering
\includegraphics[scale=0.27]{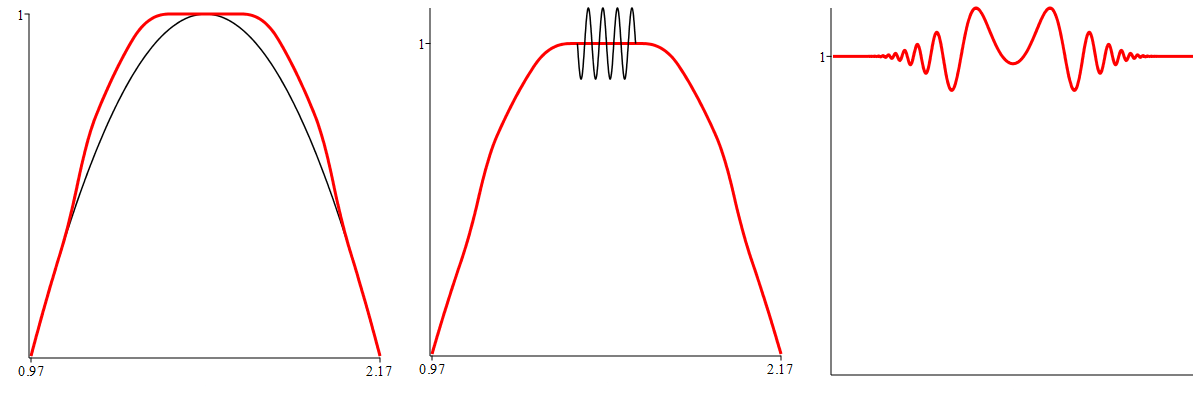}\vspace{-5mm}
\caption{Construction of a function $T_{m,n,\alpha}$}\label{F-sph-2}
\end{figure}

Figure~\ref{F-sph-2} illustrates the construction of the functions $T_{m,n,\alpha}$.
\begin{itemize}
  \item The left figure displays the graphs of the functions $\sin^m(\theta)$ (black) and $\sin^m(\theta)+P_{m,n,\alpha}(\theta)$ (red), for $\theta \in [\frac{\pi}{2}-0.6,\frac{\pi}{2}+0.6]$.
  \item The middle figure indicates (in black) where we will insert the rapidly oscillating perturbation $u_{m,n,\alpha}$ constructed in the next paragraph.
  \item The right figure displays a zoom on the function $u_{m,n,\alpha}$, whose support is contained in $\set{\theta~|~\sin^m(\theta)+P_{m,n,\alpha}(\theta)=1}$.
\end{itemize}

\FloatBarrier

\subsubsection{Construction of $u_{m,n,\alpha}$}

Define the function $v$ as follows.
\begin{equation}\label{E-S2m-16}
v(t) = \left\{
\begin{array}{l}
\exp\left( \frac{1}{(t+1)(t-1)}\right) \, \cos\left( \frac{1}{(t+1)(t-1)}\right) \text{~~for~~} -1 < t <1 \,,\\[5pt]
0 \text{~~elsewhere.}
\end{array}
\right.
\end{equation}

This function is smooth, even, bounded, with bounded first and second derivatives.
Define the family of functions $u_{m,n,\alpha} : [0,\pi] \to \R$, such that they are symmetric with respect to $\frac{\pi}{2}$, $u_{m,n,\alpha}(\pi-\theta) = u_{m,n,\alpha}(\theta)$ for all $\theta\in [0,\pi]$, and given by the formula
\begin{equation}\label{E-S2m-20}
u_{m,n,\alpha}(\frac{\pi}{2} + \hat\theta) = a_{m,n,\alpha}\, v\left( \frac{(mn)^2\, \hat{\theta}}{\alpha}\right)\,,
\end{equation}
where $\hat{\theta} \in [-\frac{\pi}{2},\frac{\pi}{2}]$. The graph of $u_{m,n,\alpha}$ appears in Figure~\ref{F-sph-2} (right). Note that $u_{m,n,\alpha}$ is supported in the set $\set{\theta ~|~ \sin^m(\theta)+P_{m,n,\alpha}(\theta) = 1}$. The constant $a_{m,n,\alpha}$ is chosen such that, for any $m,n\geq 1$, and $\alpha \in (0,\frac{1}{4}]$,
\begin{equation}\label{E-S2m-22}
|u_{m,n,\alpha}| + |u_{m,n,\alpha}'| + |u_{m,n,\alpha}''| \le \alpha \,.
\end{equation}

\subsubsection{Properties of $T_{m,n,\alpha}$}~

From the construction, $T_{m,n,\alpha}(\theta)-1$ changes sign infinitely many times on the interval $[\frac{\pi}{2} - \frac{\alpha}{(mn)^2}, \frac{\pi}{2} + \frac{\alpha}{(mn)^2}]$. Indeed, $\sin^m(\theta) + P_{m,n,\alpha}(\theta) = 1$ on that interval, and $u_{m,n,\alpha}$ changes sign infinitely often on the same interval. Also, since $\sin^m$, $P_{m,n,\alpha}$ and $u_{m,n,\alpha}$ are all smooth, $T_{m,n,\alpha}$ is smooth.

\subsection{Non-degeneracy of the metric}\label{SS-ndm2m}~

We use Subsection~\ref{SS-pe}. To the function $T_{m,n,\alpha}$ we associate the function $Q_{m,n,\alpha}$ through the relation~\eqref{E-QQm}.
This function defines a conformal metric $g_{m,n,\alpha}$ on $\Ss^2$ provided that it is positive and smooth. Taking into account the relations~\eqref{E-SHm} and \eqref{E-QQmm}. Write,
\begin{equation}\label{E-S2m-12}
1- Q_{m,n,\alpha}(\theta) = \dfrac{N(\theta)}{D(\theta)}\,,
\end{equation}
where
\begin{equation}\label{E-S2m-14}
\left\{
\begin{array}{ll}
N(\theta) &= (\cK_m P_{m,n,\alpha})(\theta) + (\cK_m u_{m,n,\alpha})(\theta)\,,\\[5pt]
D(\theta) & = m(m+1) \, \sin^2(\theta) \, T_{m,n,\alpha}(\theta)\,,
\end{array}%
\right.
\end{equation}
and recall that $\cK_m(\sin^m) \equiv 0$. Because $P_{m,n,\alpha}$ and $u_{m,n,\alpha}$ are supported in $J_{mn} := [\frac{\pi}{2}-\frac{1}{mn},\frac{\pi}{2}+\frac{1}{mn}]$, we have $Q_{m,n,\alpha} \equiv 1$ in $(0,\pi)\sm J_{mn}$.
 It therefore suffices to study $Q_{m,n,\alpha}$ in the interval $J_{mn}$ and, by symmetry with respect to $\frac{\pi}{2}$, in   $J_{+,mn} := [\frac{\pi}{2},\frac{\pi}{2}+\frac{1}{mn}]$. As above, we set $\hat\theta = \theta - \frac{\pi}{2}$.\medskip

From the definition of $T_{m,n,\alpha}$, Equation \eqref{E-S2m-8}, we deduce
\begin{equation*}
\left| T_{m,n,\alpha}(\theta) - 1\right| \le 1 - \sin^m(\theta) + \left|P_{m,n,\alpha}(\theta)\right| + \left|u_{m,n,\alpha}(\theta)\right|\,.
\end{equation*}
Using \eqref{E-S2m-22} and Proposition~\ref{P-S2-Pmn}(4), we obtain
\begin{equation*}
\left| T_{m,n,\alpha}(\theta) - 1\right| \le \frac{3m}{(mn)^2} + \alpha\,, \text{~for~} n \ge N(m,\alpha),~ \theta\in J_{+,mn}\,.
\end{equation*}

It follows that for $0 < \alpha \le \frac 14$, and $n \ge N(m,\alpha)$,
\begin{equation}\label{E-S2m-34}
  \left| T_{m,n,\alpha}(\theta) - 1\right| \le \frac{1}{2} \text{~in~} J_{+,mn}\,.
\end{equation}
In particular, this inequality implies that for $0 < \alpha \le \frac 14$,  there exists some $C_\alpha >0$ such that, for $\theta \in J_{+,mn}$ and for $n$ large enough, $D(\theta) > C_\alpha$. It follows that $Q_{m,n,\alpha}$ is well-defined on $\Ss^2$ and $C^{\infty}$ by equation \eqref{E-S2m-12}.\medskip

\textbf{Claim.~}~For $0 < \alpha \le \frac 14$, small enough, and for $n \ge N(m,\alpha)$ large enough, the function $Q_{m,n,\alpha}$ is close to $1$.\smallskip

\emph{Proof of the Claim.~}~Since $Q_{m,n,\alpha} \equiv 1$ in $(0,\pi)\sm J_{mn}$, and by symmetry around $\frac{\pi}{2}$, it suffices to consider $\theta \in J_{+,mn}$.  In this interval, we have
\begin{equation}\label{E-S2m-34a}
\big|m(m+1) \sin^2(\theta)-m^2)\big| \le m + \frac{1}{n^2}\,.
\end{equation}
From \eqref{E-S2m-34}, and for $n\ge N(m,\alpha)$ and  $\theta \in J_{+,mn}$, we obtain
\begin{equation}\label{E-S2m-34b}
D(\theta) \ge m(m+1)\, \left( 1 - \frac{5}{n^2} - \alpha \right) \ge m(m+1)(1-2\alpha)\,.
\end{equation}
We estimate $N(\theta)$, using \eqref{E-S2m-14}. From \eqref{E-Km}, \eqref{E-S2m-22},  and \eqref{E-S2m-34a}, we obtain
\begin{equation*}
\begin{array}{ll}
|(\cK_mu_{m,n,\alpha})(\theta)| & \le |u''_{m,n,\alpha}(\theta)| + |u'_{m,n,\alpha}(\theta)| \\[5pt]
& ~~~\quad + \big( m(m+1) \sin^2(\theta)-m^2)|u_{m,n,\alpha}(\theta)|\big)\\[5pt]
& \leq ( m+1) \alpha \,.
\end{array}
\end{equation*}

We estimate $|\cK_mP_{m,n,\alpha}|$ as follows.
\begin{equation*}
\begin{array}{ll}
|(\cK_mP_{m,n,\alpha})(\theta)|  & \le |P''_{m,n,\alpha}(\theta)| + |P'_{m,n,\alpha}(\theta)|\\[5pt] &~~~+ \big( m(m+1) \sin^2(\theta)-m^2)|P_{m,n,\alpha}(\theta)|\big) \,.
\end{array}%
\end{equation*}

Using the estimates in Proposition~\ref{P-S2-Pmn}, and the fact that $|\sin(t)| \le |t|$, we obtain the following inequalities for $n \ge N(m,\alpha)$ and $\theta \in J_{+,mn}$,
\begin{equation*}
\left\{
\begin{array}{rl}
|P''_{m,n,\alpha}(\theta)| &\le m(1+5\alpha)\,,\\[5pt]
|P'_{m,n,\alpha}(\theta)| &\le \frac{2m}{(mn)^2}\,,\\[5pt]
|P_{m,n,\alpha}(\theta)| &\le \frac{2m}{(mn)^3}\,.
\end{array}%
\right.
\end{equation*}

From these estimates and \eqref{E-S2m-34a}, we obtain
\begin{equation*}
|(\cK_mP_{m,n,\alpha})(\theta)| \le m(1+5\alpha) + \frac{5}{n^2}
\end{equation*}
for $n \ge N(m,\alpha)$ and $\theta \in J_{+,mn}$.

\medskip

Finally, for $\theta \in J_{+,mn}$ and $n \ge N(m,\alpha)$, we have
\begin{equation*}
\left\{
\begin{array}{rl}
 m(m+1)(1-2\alpha) &\le D(\theta)\,,\\[5pt]
N(\theta) &\le  m \, \left( 1 + \frac{5}{n^2} + 5 \alpha \right)\,,\\[5pt]
\big| 1 - Q_{m,n,\alpha}(\theta)| &\le \frac{1}{m+1} \, \left( 1 + \frac{10}{n^2} + 10 \alpha \right)\,.
\end{array}%
\right.
\end{equation*}

The claim is proved.\hfill\qed\medskip

Combining the previous estimates, we obtain the main result of this section.

\begin{proposition}\label{P-S2m-6}
For any $m\ge 1$, and $\alpha \in (0,\frac{1}{24})$, there exists $N_1(m,\alpha)$ such that, for $n\geq N_1(m,\alpha)$,
\begin{equation}\label{E-Q6}
| Q_{m,n,\alpha} (\theta)-1| \leq \frac{1}{m+1}(1+12\alpha)\,.
\end{equation}
In particular, the metric $g_{m,n,\alpha} = Q_{m,n,\alpha}\, g_0$ is smooth, non-degenerate, and close to $g_0$ provided that $m$ is large enough.
\end{proposition}

\subsection{Proof of Propositions~\ref{P-S2m-2} and \ref{P-S2m-8}}

We  now apply the results obtained in Subsections~\ref{SS-pe}, \ref{SS-psm}, and \ref{SS-ndm2m}.

\subsubsection{Proof of Proposition~\ref{P-S2m-2}} Fix $m=1$, and define the function $\Phi$ in spherical coordinates by,
\begin{equation*}
\Phi(\theta,\phi) = T_{1,n,\alpha}(\theta) \, \cos(\phi)\,,
\end{equation*}
where $0 < \alpha \le \frac{1}{24}$, and $n$ is large enough, according to Subsection~\ref{SS-ndm2m}. The function $\Phi$ is clearly smooth away from the north and south poles of the sphere ($\theta$ away from $0$ and $\pi$). Near the poles, $\Phi$ is equal to the spherical harmonic $\sin(\theta)\,\cos(\phi)$. It follows that $\Phi$ is smooth. By Proposition~\ref{P-S2m-6}, the function $Q_{1,n,\alpha}$ associated with $T_{1,n,\alpha}$ by the relation \eqref{E-QQm} extends to a smooth positive function on $\Ss^2$. Choose $G = Q_{1,n,\alpha}\, g_0$, where $g_0$ is the standard round metric. Then, according to Subsection~\ref{SS-pe}
\begin{equation}\label{E-S2m-50}
-\Delta_{G} \Phi = 2\, \Phi.
\end{equation}
This proves Assertions~\ref{P-S2m-2}(2) and (3). Assertion~\ref{P-S2m-2}(4) is a consequence of the min-max. Indeed, by Proposition~\ref{P-S2m-6},
\begin{equation}\label{E-S2m-50b}
\frac 12 - 6 \alpha \le G \le \frac 32 + 6\alpha\,.
\end{equation}
According to our choice of $\alpha$, the left-hand side of this inequality is positive. Call $R_G$, resp. $R_0$, the Rayleigh quotient of $(\bS^2,g_G)$, resp. on $(\bS^2,g_0)$. Then, by (\ref{E-S2m-50b}),
\begin{equation*}
\frac 23\,(1 + 4\alpha)^{-1} \, R_0(\psi) \le R_G(\psi) \le \frac 12\, (1 - 12 \alpha)^{-1} \, R_0(\psi)\,,
\end{equation*}
for all $0 \not = \psi \in \bS^2$. From the min-max, we conclude that
\begin{equation}\label{E-S2m-52}
\frac 23 \, (1 + 4\alpha)^{-1} \, \lambda_k(g_0) \le \lambda_k(g_G) \le 2\, (1 - 12 \alpha)^{-1}\, \lambda_k(g_0)\,,
\end{equation}
for all $k\ge 1$, where $\lambda_k(g)$ denotes the $k$-th eigenvalue of the Laplace-Beltrami operator for the metric $g$ (eigenvalues arranged in nondecreasing order, starting from the labelling $1$, with multiplicities accounted for).\medskip

We have $\lambda_1(g_0) = \lambda_1(g_G) = 0$. Since $\lambda_2(g_0) = \cdots = \lambda_4(g_0) = 2$, and $\lambda_5(g_0) = \cdots = \lambda_9(g_0) = 6$, we conclude from \eqref{E-S2m-52} and our choice of $\alpha$, that
\begin{equation}\label{E-S2m-54}
 2 < 4\, (1 + 4\alpha)^{-1}  \le \lambda_5(g_G)\,.
\end{equation}
From \eqref{E-S2m-50} we know that the eigenfunction $\Phi$ is associated with the eigenvalue $2$ of $-\Delta_G$. It follows from \eqref{E-S2m-54} that this eigenvalue has labelling at most $4$. This proves Assertion~\ref{P-S2m-2}(4). \medskip

Assertion~\ref{P-S2m-2}(1) is similar to Assertion~(1) in Proposition~\ref{P-S2m-8}. We defer its proof to Paragraph~\ref{PP}. ~~Proposition~\ref{P-S2m-2} is proved. \hfill \qed \medskip

\subsubsection{Proof of Proposition~\ref{P-S2m-8}} According to Subsection~\ref{SS-ndm2m}, when $m\ge 1$, an appropriate choice of $(\alpha,n)$ yields a function
$$
\Phi_m(\theta,\phi) = T_{m,n,\alpha}(\theta) \, \cos(m\phi)\,,
$$
and a function $G_m = Q_{m,n,\alpha}$ satisfying \eqref{E-Q6}, such that
\begin{equation*}
-\Delta_{G_m} \Phi_m = m(m+1) \Phi_m.
\end{equation*}
Choosing $m$ large enough, the metric $g_m = G_m\, g_0$ can be made as close as desired to the standard metric $g_0$, see however Remark~\ref{R-cg0}. This proves Assertions~\ref{P-S2m-8}(2)--(4). \medskip

\subsubsection{Proof of Assertions~\ref{P-S2m-2}(1) and \ref{P-S2m-8}(1) \emph{(Nodal properties of the  eigenfunction $\Phi_m$)}}\label{PP}

~~For simplicity, denote the function $T_{m,n,\alpha}$ by $T$, so that
$$
T(\theta) = \sin^m(\theta) + P_{m,n,\alpha}(\theta) + u_{m,n,\alpha}(\theta).
$$
Let $V$ denote the function
$$
V(\theta) = \sin^m(\theta) + P_{m,n,\alpha}(\theta).
$$

Taking Proposition~\ref{P-S2-Pmn} into account, we have the following properties for $V$: $V(\pi - \theta) = V(\theta)$, $V(\theta) = \sin^m(\theta)$ near $0$ and $\pi$, and $V \equiv 1$ for $|\theta-\frac{\pi}{2}|\le \frac{\alpha}{(mn)^2}$. In the interval $[\frac{\pi}{2},\pi]$, we have,
$$
V'(\theta) = m\, \cos(\theta)\,\sin^{m-1}(\theta)\,\left( 1 - \chi_{\alpha}(m^2n^2\hat\theta)\right) -\, \gamma\; \xi(mn\, \hat\theta)\,,
$$
where $\gamma := \gamma_{m,n,\alpha}$, see \eqref{E-S2m-gam}.
It follows that $V'(\theta) \le 0$ in $[\frac{\pi}{2},\pi]$, so that $0\le V(\theta) \le 1$ in $[0,\pi]$. With the notation $u=u_{m,n,\alpha}$, recalling that $V(\theta) \ge 0$, and that $\supp(u) \subset \set{V=1}$, we conclude that
\begin{equation}
\left\{
\begin{array}{l}
\set{u(\theta) < 0}\times [0,2\pi] \subset \set{\Phi_m < 1}\,,\\[5pt]
\set{\Phi_m \ge 1} \subset \set{u(\theta) \ge 0}\times [0,2\pi]\,,\\[5pt]
\set{u(\theta) > 0}\times \set{0} \subset \set{\Phi_m > 1}\,.
\end{array}%
\right.
\end{equation}

This means that the set $\set{\Phi_m > 1}$ has at least one connected component in each band $\set{u(\theta) > 0}\times [0,2\pi]$. \smallskip

The proof of Proposition~\ref{P-S2m-8} is now complete.\hfill \qed

\begin{remark}\label{R-cg0}
Note that, by \eqref{E-S2m-P1}, for any $n$ and $\alpha$, $T_{m,n,\alpha}(\pi/2)=1$, $T'_{m,n,\alpha}(\pi/2)=0$, $T''_{m,n,\alpha}(\pi/2)=0$, and hence, by the relation \eqref{E-QQm}, $Q_{m,n,\alpha}(\pi/2) = \frac{m}{m+1}\neq 1$. Therefore, it is impossible for each $G_m$ to be arbitrarily close to the round metric, regardless of the choice of $n$ and $\alpha$. Proposition \ref{P-S2m-8} merely states that we can find a sequence of metrics that converge in some sense to the round metric.
\end{remark}

\begin{remark}
If we restrict the eigenfunction $\Phi$ (Proposition~\ref{P-S2m-2}) or $\Phi_m$ (Proposition~\ref{P-S2m-8}) to the hemisphere $\set{(\theta,\phi) \in [0,\frac{\pi}{2})\times [0,2\pi]}$, we obtain counterexamples to Conjecture~\ref{C-ECP}  for the hemisphere, with a metric conformal to the standard metric $g_0$, and Neumann boundary condition.
\end{remark}%

\section{Final comments}\label{S-fin}

With respect to the Extended Courant property, we would like to point out that there are ways of counting nodal domains of sums of eigenfunctions which avoid the pathologies exhibited in the examples constructed in Sections~\ref{S-T2} and \ref{S-S2m}. In the deterministic framework, we mention \cite{PoSo2007,PoPoSt2019} in which the nodal count involves some weights. In the probabilistic framework, the topological complexity of the nodal set of a random sum of eigenfunctions can be estimated. We refer to the recent thesis \cite{Riv2018}, and its bibliography.\medskip

With respect to \cite{BLS}, we would like to point out that although our starting point is the same (the idea to construct Liouville metrics with an oscillatory component), our goals and methods are different.

\vspace{1cm}
\appendix

\section{Bounds on the number of nodal domains on $\Ss^2$ with the round metric}\label{A-S2}

The following result can be found in \cite{Char2015}:

\begin{proposition}\label{boundsphere}
Let $ f: \R^n \to \R$ be a polynomial of degree $d$. Then, the number of nodal domains of its restriction to $\Ss^{n-1}$ is bounded by $2^{2n-1}d^{n-1}$.
\end{proposition}

In the case of $\mathbb{S}^2$ with the round metric, every eigenfunction is  the restriction of a harmonic homogeneous polynomial to the sphere. Also, for such a polynomial of degree $ \ell$, its eigenvalue on the sphere is $\ell(\ell+1)$,  with multiplicity $2\ell +1$. For a sum $w$ of spherical harmonics of degree less than or equal to $\ell$, Conjecture~\ref{C-ECP} would give $$\beta_0(w) \leq 1 + \sum\limits_{k=0}^{l-1}(2k+1) \le \ell^2+1.$$ Using Proposition \ref{boundsphere},  we get the following weaker estimate.

\begin{corollary}
Let $g_0$ be the round metric for $\Ss^2$. Then, the sum $w$ of spherical harmonics of degree less than or equal to $\ell$ has at most $8\, \ell^2$ nodal domains.
\end{corollary}

However, the direct nodal count is highly unstable in the case of $C^\infty$ metrics, as we have shown in  Section~\ref{S-S2m},  see also Section~\ref{S-fin}.

\section{Isotropic quantum harmonic oscillator in dimension $2$}\label{A-QHO}

In this section, we will show that Conjecture \ref{C-ECP} is true for the harmonic oscillator $H : L^2(\mathbb{R}^2) \to L^2(\mathbb{R}^2), H = -\Delta + x^2 + y^2$.

\begin{proposition}\label{ECPQHO}
	Let $f_i$ be the eigenfunctions of $H$ with eigenvalues ordered in increasing order with multiplicities. Then, for any linear combination $f = \sum\limits_{i=1}^{n}a_i f_i$, we have $ \beta_0(f) \leq n\,$.
\end{proposition}

A basis $f_n$ of eigenfunctions of $H$ is given by  $$
H_{a,b}(x,y):= e^{-\frac{x^2 + y^2}{2}}H_{a}(x)H_b(y)\,, \,0 \leq a,b \in \N\,,$$
 where $H_n$ refers to the $n$-th Hermite polynomial.\\
  The associated eigenvalue is given by $2(a+b+1)$,  with multiplicity $a+b+1$. Therefore, counting multiplicities, for each $n$ in the interval $[\frac{k(k+1)}{2} +1, \frac{(k+1)(k+2)}{2}]$ for some positive integer $k$, $f_n$ is a polynomial of degree $k$.\\

 For a polynomial $f$ of degree $k$ in $2$ variables, we have the following upper bound on the number of its nodal domains:

\begin{lemma}\label{upperboundpoly}
	For any polynomial $f$ of degree $k$ in $\mathbb{R}^2$, $$ \beta_0(f)\leq k(k+1)/2 +1\,.$$
	The upper bound is achieved by $k$ non-parallel lines.
\end{lemma}

To prove this, we first note that the number of nodal domains is bounded from above by $U(f)+ S(f)+1$, where $U(f)$ is the number of connected components of the nodal set and
 $ S(f)= \sum (s_i-1)$, where the sum is taken over all singular points $a_i$ and $s_i$ is the order of the singularity at $a_i$ (the lowest homogeneous order term in the Taylor expansion of $f$ around $a_i$ ).

Now, we use classical theorems by B\'{e}zout and Harnack, see \cite{BCR}. Recall that for a curve $\gamma$ defined by $\gamma = F^{-1}(0)$ for some polynomial $F$, a singular point of degree $d$ is a point $x$ on $\gamma$ such that all partial derivatives of $F$ of order less than or equal to $d$ vanish at $x$, but some derivative of order $d+1$ does not vanish.

\begin{theorem}[B\'{e}zout's theorem]
Let $f$ and $g$ be real algebraic curves of degree $m$ and $n$. If the number of points in the intersection of $f$ and $g$ is infinite, then the polynomials defining $f$ and $g$ have a common divisor. If the number of points in the intersection of $f$ and $g$ is finite, then it is less than or equal to $mn$.
\end{theorem}

\begin{theorem}[Harnack's curve theorem]
 Let $f$ be a real irreducible polynomial in two variables, of degree $k$. Let $a_i$ be the singular points of the nodal set, with order $s_i$. We have the following inequality\footnote{In fact, the original theorem as stated in \cite{BCR} deals with algebraic curves in  $\mathbb{RP}^2$. However, it is easily adapted to $\R^2$ by adding at most $k$ unbounded components. } for the number of connected components of its nodal set:

 $$U(f) \leq \frac{(k-1)k}{2} - \sum_i \frac{s_i(s_i-1)}{2} + 1\,.$$

\end{theorem}

Now, we proceed by induction. For $k=1$, the lemma is trivial. Now, consider a polynomial $f$ of degree $k >1$. It can be either irreducible or the product of two smaller degree polynomials.

If $f$ is irreducible, then by Harnack's theorem we have $$\beta_0(f) \leq (k-1)k/2 +2 \,,$$ since for all $a \geq 1$, $a-1 \leq a(a-1)/2\,$.

If $f = PQ$ with $\deg P = j$ and $\deg Q = k-j$, the number of nodal domains is bounded by $\beta_0(P)+\beta_0(Q) + j(n-j) -1\,$. Indeed, every intersection between $P$ and $Q$ adds the same number of nodal domains as the degree of their intersection, and this number can be  bounded  by B\'{e}zout's theorem. We need to substract $1$ to remove the initial original domain of $\mathbb{R}^2$ (otherwise, multiplying two linear functions would give $5$ nodal domains.)

By induction, we have the following inequality:

\begin{align*}
	\beta_0(f) &\leq \frac{j(j+1)}{2} + \frac{(n-j)(n-j+1)}{2} + j(n-j) + 1 \\
	&\leq \frac{n(n+1)}{2}+1 \,.
\end{align*}

Now, since this was achieved by $P$ and $Q$ being the product of linear factors, then $f$ is a product of linear factors. This proves lemma \ref{upperboundpoly}.\hfill\qed\medskip

We can now complete the proof of proposition \ref{ECPQHO}.\\  Let $n \in [k(k+1)/2 +1, (k+1)(k+2)/2]\,$. Then, any linear combination of $f_1, f_2, \ldots, f_n$ will be a polynomial of degree at most $k$. Any such polynomial has at most $k(k+1)/2 + 1$ nodal domains. Therefore, Conjecture \ref{C-ECP} is true in the case of the isotropic two-dimensional quantum harmonic oscillator.

\begin{remark}
	It is still unclear if this upper bound can be reached for any $k >2\,$.
\end{remark}
\begin{remark}
	 Considering the results of this paper, it seems likely that a small perturbation of either the metric in $\mathbb{R}^2$ or the potential could break this upper bound.
\end{remark}

\vspace{1cm}
\bibliographystyle{plain}

\vspace{2cm}
\end{document}